\documentclass{mod_elsart}

\usepackage{amsmath}
\usepackage{amssymb}
\usepackage{graphics}
\usepackage{graphicx}
\usepackage{epsfig}
\usepackage{setspace}
\usepackage{amssymb}

\begin{document}

\begin{frontmatter}



\title{Sulfate attack  in sewer pipes: Derivation of a concrete corrosion model via two-scale convergence}
\vspace{-0.3in}
\author[label1]{Tasnim Fatima} and \author[label2]{Adrian Muntean}

\address[label1]{Centre for Analysis, Scientific computing and Applications (CASA),\\
Department of Mathematics and Computer Science,\\
Technical University Eindhoven, Eindhoven, The Netherlands.\\
E-mail: t.fatima@tue.nl}

\address[label2]{Centre for Analysis, Scientific computing and Applications (CASA),\\
Institute for Complex Molecular Systems (ICMS),\\
Department of Mathematics and Computer Science,\\
Technical University Eindhoven, Eindhoven, The Netherlands.\\
E-mail: a.muntean@tue.nl}

\begin{abstract}
We explore the homogenization limit and rigorously derive upscaled equations for a microscopic
 reaction-diffusion system modeling sulfate corrosion in sewer pipes made of concrete.
The system, defined in a periodically-perforated domain, is
semi-linear, partially dissipative and weakly coupled via a
non-linear ordinary differential equation posed on the solid-water interface at the pore level.
Firstly, we show the well-posedness of the microscopic model. We
then apply homogenization techniques based on two-scale convergence
for an uniformly periodic domain and derive upscaled equations
together with explicit formulae for the effective diffusion
coefficients and reaction constants. We use a boundary unfolding
method to pass to the homogenization limit in the non-linear ordinary differential equation. Finally, besides giving its strong formulation, we also
prove that the upscaled two-scale model admits a unique solution.
\end{abstract}
\vspace{-0.3in}
\end{frontmatter}
\vspace{-0.5in}
\section{Introduction}
\vspace{-0.3in}

This paper treats the periodic homogenization of a semi-linear reaction-diffusion system coupled with a nonlinear differential equation arising in
 the modeling of the sulfuric acid attack in sewer pipes made of concrete. The concrete corrosion situation we are dealing with here strongly influences the durability of
 cement-based materials especially in hot environments leading to spalling of concrete and macroscopic fractures of sewer pipes. It is financially important to have
 a good estimate on the moment in time when such pipe systems need to be replaced, for instance, at the level of a city like Los Angeles.
 To get good such practical estimates, one needs
 on one side easy-to-use macroscopic corrosion models to be used for a numerical forecast of corrosion, while on the other side one needs to ensure
 the reliability of the averaged models by allowing them to incorporate a certain amount of microstructure information. The relevant question is: {\em How much of this oscillatory-type information is needed to get a sufficiently accurate description of the heterogeneous reality?} Due to the complexity of
 possible shapes of the microstructure, averaging concrete materials is far more difficult than averaging metallic composites with rigorously defined
 well-packed structure.  In this paper, we imagine our concrete piece to be made of a periodically-distributed microstructure. Based on this assumption, we
 provide here a rigorous justification of the formal asymptotic expansion performed by us (in
  \cite{tasnim1}) for this reaction-diffusion scenario. Note that in  \cite{tasnim1}  upscaled models are derived for a more general situation involving a locally-periodic distribution of perforations\footnote{The word "perforation" is seen here as a synonym for  "pore"
  or "microstructure".}. Locally periodic geometries refer to a special case of $x$-dependent microstructures, where, inherently, the outer normals to
(microscopic) inner interfaces are dependent on both spatial slow variable, say $x$, and fast variable, say $y$.

In the framework of this paper, we combine two-scale convergence
concepts with the periodic unfolding of interfaces to pass to the
homogenization limit (i.e. to $\varepsilon\to 0$,
 where $\varepsilon$ is a small parameter linked to the relative size of the perforation) for the uniformly periodic case. Here, the outer normals
  to the inner interfaces are dependent only on the spatial fast variable. For more details on the mathematical modeling of sulfate corrosion
  of concrete, we refer the reader to \cite{mbp,Niko} (a moving-boundary  approach: numerics and formal matched asymptotics), \cite{tasnim2}
  (a two-scale reaction-diffusion system modeling sulfate corrosion), as well as to \cite{maria1}, where a nonlinear Henry-law type transmission condition
   (modeling $H_2S$ transfer across all air-water interfaces present in this sulfatation problem) is analyzed. Mathematical background on  periodic
    homogenization  can be found in e.g., \cite{Papanicolau,Cioranescu,Persson},  while a few relevant (remotely resembling) worked-out examples
     of this averaging  methodology are explained, for instance, in
\cite{Hornung2,Belyaev,Peter,Meier,Meier1,Noorden}. It is worth
noting that, since it deals with the homogenization of a linear
Henry-law setting, the paper \cite{Peter} is related to our
approach. The major novelty here compared to  \cite{Peter} is that
we now need to pass to the limit in a non-dissipative object, namely
a nonlinear ordinary differential equation (ode). The ode is
describing sulfatation reaction at the inner water-solid interface
-- place where corrosion localizes. This aspect  makes a rigorous
averaging challenging. For instance, compactness-type methods do not
work in the case when the nonlinear ode is posed on
$\epsilon$-dependent surfaces. We circumvent this issue by "boundary
unfolding" the ode. Thus we fix,  as independent of $\epsilon$, the
reaction interface similarly as in \cite{MM}, and only then we pass
to the limit. Alternatively, one could use varifolds (cf. e.g.
\cite{varifolds}), since this seems to be the natural framework for
the rigorous passage to the limit when both the surface measure and
the oscillating sequences depend on $\epsilon$. However, we find the
boundary unfolding technique easier to adapt to our scenario than
the varifolds.

Note that here we approach the corrosion problem deterministically.
However, we have reasons to expect that the uniform periodicity
assumption can be relaxed by assuming instead a Birkhoff-type
ergodicity of the microstructure shapes and positions, and hence,
the natural averaging context seems to be the one offered by random
fields; see ch. 1, sect. 6 in \cite{Chechkin2}, ch. 8 and 9 in
\cite{Jikov}, or \cite{BMP}. But, methodologically, how big is the
overlap between homogenizing deterministically locally-periodic
distributions of microstructures compared to working in the random
fields context? We will treat these and related aspects elsewhere.

The paper is organized as follows: We start off in section \ref{micro} (and continue in section \ref{basic}) with the analysis of the microscopic model.
 In section \ref{est}, we obtain the $\varepsilon$-independent estimates needed for the passage to the limit $\varepsilon\to 0$.  Section \ref{two-scale}
  contains the main result of the paper: the set of the upscaled two-scal equations.
\vspace{-0.3in}
\section{The microscopic model}\label{micro}
\vspace{-0.3in}
In this section, we  describe the geometry of our array of periodic microstructures and briefly indicate the most aggressive chemical reaction mechanism typically active in sewer pipes. Finally, we list the set of microscopic equations.
\vspace{-0.3in}
\subsection{Basic geometry}\label{geometry}
\begin{figure}[htbp]
\begin{center}
\includegraphics[width=\textwidth]{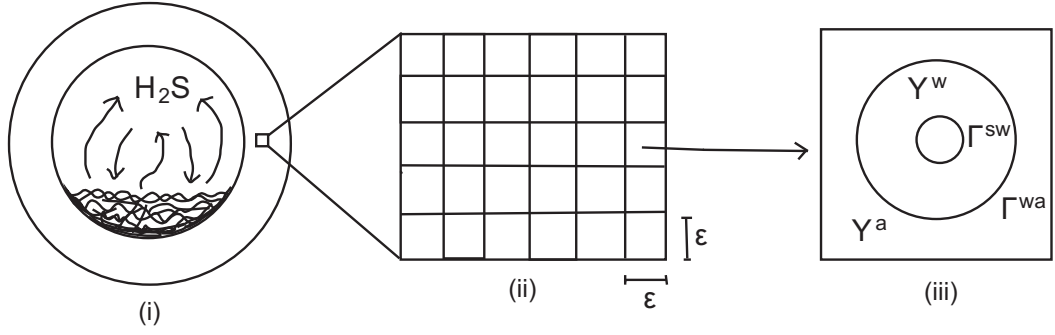}\\
\caption{\small{Left: Cross-section of a sewer pipe pointing out one region. Middle:
 Periodic approximation of the periodic rectangular domain. Right: Reference pore configuration.}}\label{fig1}
\end{center}
\end{figure}

\vspace{-0.3in}
Fig. \ref{fig1} (i) shows a cross-section of a sewer pipe hosting
corrosion. We assume that the geometry of the porous
medium in question consists of a system of pores periodically
distributed inside the three-dimensional cube $\Omega:=[a,b]^3$ with
$a,b\in \mathbb{R}$ and $b> a$.  The exterior boundary of
$\Omega$   consists of two disjoint, sufficiently smooth parts: ${
\Gamma^N}$ - the Neumann boundary and ${\Gamma^D}$ - the Dirichlet
boundary. The reference pore, say $Y:=[0,1]^3$, has three pairwise disjoint connected domains $Y^s$, $Y^w$ and $Y^a$ with smooth boundaries
$\Gamma^{sw}$ and $\Gamma^{wa}$,  as shown in Fig. \ref{fig1} (iii).
Moreover, $Y:=\bar{Y}^s\cup \bar{Y}^w\cup \bar{Y}^a$.

Let $\varepsilon$ be a sufficiently small scaling factor denoting the
ratio between the characteristic length of the pore $Y$ and the
characteristic length of the domain $\Omega$. Let
$\chi^w$ and $\chi^a$ be the characteristic functions of the sets
$Y^w$ and $Y^a$, respectively. The shifted set $Y^w_k$ is defined by
$$Y^w_k:=Y+\Sigma_{j=0}^3k_je_j\; \mbox{for}\;\;k=(k_1,k_2,k_3)\in \mathbb{Z}^3,$$
where $e_j$ is the $j^{th}$ unit vector. The union of all shifted
subsets of $Y^w_k$ multiplied by $\varepsilon$ (and confined
within $\Omega$) defines the perforated domain $\Omega^\varepsilon$,
namely
$$\Omega^\varepsilon:= \cup_{k\in \mathbb{Z}^3} \{\varepsilon Y^w_k\;|\; \varepsilon Y^w_k \subset \Omega\}.$$
Similarly, $\Omega_1^\varepsilon$, $\Gamma^{sw}_\varepsilon$, and
$\Gamma^{wa}_\varepsilon$ denote the union of the shifted subsets
(of $\Omega$) $Y^a_k$, $\Gamma^{sw}_k$, and $\Gamma^{wa}_k$ scaled
by $\varepsilon$.
Since usually the concrete in sewer pipes is not completely dry, we
decide to take into account a partially saturated porous
material\footnote{The solid, water and air parts corresponds to $Y^s$, $Y^w$ and $Y^a$,
respectively.}. We assume that every pore has three distinct
non-overlapping parts: a solid part (grain) which is placed in the
center of the pore, the water film which surrounds the solid part,
and an air layer bounding the water film and filling the space of
$Y$ as shown in Fig. \ref{fig1}. The air connects neighboring
pores to one another. The geometry defined above satisfies the
following assumptions:
\vspace{-0.15in}
\begin{enumerate}
  \item Neither solid  nor water-filled parts touch the boundary of the pore.
  \item All internal (air-water and water-solid) interfaces are sufficiently smooth and do not touch each other.
\end{enumerate}
\vspace{-0.2in}
These geometrical restrictions imply that the pores are connected by air-filled parts only which is needed not only to give a meaning
to functions defined across interfaces, but also  to introduce the
concept of extension as given, for instance, in \cite{paulin}. Furthermore, there are no solid-air interfaces.
\vspace{-0.3in}
\subsection{Description of the chemistry}
\vspace{-0.3in}
There are many variants of severe attack to concrete in sewer pipes,
we focus here on the most aggressive one -- the sulfuric acid
attack. The situation can be described briefly as follows: (The
anaerobic bacteria in the flowing waste water release hydrogen
sulfide gas ($H_2S$) within the air space of the pipe. These
bacteria are especially active in hot environments. From the air
space inside the pipe, $H_2S(g)$\footnote{$H_2S(g)$ and $H_2S(aq)$ refer  to  \underline{g}aseous, and respectively,  \underline{aq}ueous $H_2S$.} enters the pores of the concrete
matrix where it diffuses and then dissolves in the pore water. The
aerobic bacteria catalyze some of the $H_2S$ into sulfuric acid
$H_2SO_4$. $H_2S$ molecules can move between  air-filled part and water-filled part the  water-air interfaces \cite{Balls}. We model this
microscopic interfacial transfer via Henry's law
\cite{Danckwerts}, (see the boundary conditions at
$\Gamma^{wa}_\varepsilon$ in (\ref{ii}) and (\ref{iii})). $H_2SO_4$
being an aggressive acid reacts with the solid matrix\footnote{The
solid matrix is assumed here to consist of $CaCO_3$ only. This
assumption can be removed in the favor of a more complex cement
chemistry.}
 at the solid-water interface,
 which is made up of cement, sand, and aggregate, and produces gypsum (i.e. $CaSO_4\cdot
 2H_2O$). Here we restrict our attention to a minimal set of chemical reactions
 mechanisms as suggested in \cite{mbp}, namely.
\vspace{-0.1in}
\begin{equation}\label{reeq}
 \left\{\begin{aligned}
10H^{+}+SO_4^{-2}+\mbox{org. matter} \;&\longrightarrow\;H_2S(aq)+4H_2O
+\mbox{oxid. matter}\\
H_2S(aq)+2O_2\;&{\longrightarrow}\;2H^++SO_4^{-2}\\
H_2S(aq)\;&\rightleftharpoons \;H_2S(g)\\
2H_2O+H^++SO_4^{-2}+CaCO_3\;&{\longrightarrow}\;CaSO_4\cdot 2H_2O+HCO_3^-
\end{aligned}
\right.
\end{equation}
We assume that reactions (\ref{reeq}) do not interfere with the mechanics of
the solid part of the pores. This is a rather strong assumption
since it is known that (\ref{reeq}) can actually produce local
ruptures of the solid matrix \cite{taylor}. For more details on the
involved cement chemistry and connections to acid corrosion, we
 refer  the reader to \cite{beddoe1} (for a nice enumeration of the involved physicochemical
 mechanisms),
\cite{taylor} (standard textbook on cement chemistry), as well as to
\cite{beddoe2,beddoe3,tixier} and references cited therein. For a
mathematical approach of a similar theme related to the conservation
and restoration of historical monuments, we refer to the work by R.
Natalini and co-workers (cf. e.g. \cite{natalini}).
\vspace{-0.3in}
\subsection{Setting of the equations}
\vspace{-0.3in}

The data and unknown are given by
\begin{equation*}
\begin{aligned}
{u^\varepsilon_1}_0& : \Omega \longrightarrow \mathbb{R}_+\mbox{ - initial concentration of }H_2SO_4(aq)\\
{u^\varepsilon_2}_0& : \Omega \longrightarrow\mathbb{R}_+\mbox{ -  initial concentration of }H_2S(aq)\\
{u^\varepsilon_3}_0 &: \Omega \longrightarrow \mathbb{R}_+\mbox{ - initial concentration of }H_2S(g)\\
{u^\varepsilon_4}_0 &: \Omega \longrightarrow \mathbb{R}_+\mbox{ - initial concentration of moisture}\\
{u^\varepsilon_5}_0 &: \Omega \longrightarrow \mathbb{R}_+\mbox{ - initial concentration of gypsum}\\
u_3^D&:\Gamma_D\times (0,T) \longrightarrow \mathbb{R}_+\mbox{ - exterior concentration (Dirichlet data) of }H_2S(g)\\
u^\varepsilon_1&:\Omega^\varepsilon\times (0,T)\longrightarrow \mathbb{R}\mbox{ - concentration of }H_2SO_4(aq)\\
u^\varepsilon_2&:\Omega^\varepsilon\times (0,T)\longrightarrow \mathbb{R}\mbox{ - concentration of }H_2S(aq)\\
u^\varepsilon_3&:\Omega^\varepsilon_1\times (0,T)\longrightarrow \mathbb{R}\mbox{ - concentration of }H_2S(g)\\
u^\varepsilon_4&:\Omega^\varepsilon\times (0,T)\longrightarrow \mathbb{R}\mbox{ - concentration of moisture} \\
u^\varepsilon_5&:\Gamma^{sw}_\varepsilon\times (0,T)\longrightarrow \mathbb{R}\mbox{ - concentration of gypsum } \end{aligned}
\end{equation*}
All concentrations are viewed as mass concentrations.
We consider the following system of mass-balance equations defined at the pore level.
The mass-balance equation for $H_2SO_4$ is
\begin{equation}\label{i}
\begin{aligned}
\partial_t u^\varepsilon_1+div(-d^\varepsilon_1\nabla u^\varepsilon_1)&=-k^\varepsilon_1u^\varepsilon_1+
k^\varepsilon_2u^\varepsilon_2,\quad x\in \Omega^\varepsilon, \;t \in (0,T)\\
 u^\varepsilon_1(x,0)&={u^\varepsilon_1}_0(x),\quad x\in \Omega^\varepsilon\\
-n^\varepsilon \cdot d^\varepsilon_1\nabla u^\varepsilon_1&=0,\quad x\in \Gamma^{wa}_\varepsilon,\;t \in (0,T)\\
-n^\varepsilon\cdot d^\varepsilon_1\nabla u^\varepsilon_1&=\varepsilon\eta(u^\varepsilon_1,u^\varepsilon_5),
\quad x\in \Gamma^{sw}_\varepsilon\;t \in (0,T).
 \end{aligned}
\end{equation}
The mass-balance equation for $H_2S(aq)$ is given by
\begin{equation}\label{ii}
\begin{aligned}
\partial_t u^\varepsilon_2+div(-d^\varepsilon_2\nabla u^\varepsilon_2)&=k^\varepsilon_1u^\varepsilon_1-
k^\varepsilon_2u^\varepsilon_2,\quad x\in \Omega^\varepsilon,\;t \in (0,T),\\
 u^\varepsilon_2(x,0)&={u^\varepsilon_2}_0(x),\quad x\in \Omega^\varepsilon\\
-n^\varepsilon\cdot d^\varepsilon_2\nabla u^\varepsilon_2&=\varepsilon (a^\varepsilon (x)u^\varepsilon_3-
 b^\varepsilon (x)u^\varepsilon_2),\quad x\in \Gamma^{wa}_\varepsilon,\;t \in (0,T)\\
-n^\varepsilon\cdot d^\varepsilon_2\nabla u^\varepsilon_2&=0,\quad x\in \Gamma^{sw}_\varepsilon,\;t \in (0,T).
 \end{aligned}
\end{equation}
The mass-balance equation for $H_2S(g)$ reads
\begin{equation}\label{iii}
\begin{aligned}
\partial_t u^\varepsilon_3+div(-d^\varepsilon_3\nabla u^\varepsilon_3)&=
0,\quad x\in \Omega^\varepsilon_1,\;t \in (0,T)\\
 u^\varepsilon_3(x,0)&={u^\varepsilon_3}_0(x),\quad x\in \Omega^\varepsilon_1\\
-n^\varepsilon\cdot d^\varepsilon_3\nabla u^\varepsilon_3&=0,\quad x\in \Gamma^{N},\;t \in (0,T)\\
 u^\varepsilon_3(x,t)&= u^D_3(x,t),\quad x\in \Gamma^{D},\;t \in (0,T)\\
-n^\varepsilon\cdot d^\varepsilon_3\nabla u^\varepsilon_3&=-\varepsilon (a^\varepsilon (x)u^\varepsilon_3
-b^\varepsilon (x)u^\varepsilon_2),\quad x\in \Gamma^{wa}_\varepsilon,\;t \in (0,T).
 \end{aligned}
\end{equation}
The mass-balance equation for moisture follows
\begin{equation}\label{iv}
\begin{aligned}
\partial_t u^\varepsilon_4+div(-d^\varepsilon_4\nabla u^\varepsilon_4)&=k^\varepsilon_1u^\varepsilon_1,\quad x\in \Omega^\varepsilon,\;t \in (0,T) \\
 u^\varepsilon_4(x,0)&={u^\varepsilon_4}_0(x),\quad x\in \Omega^\varepsilon\\
-n^\varepsilon\cdot d^\varepsilon_4\nabla u^\varepsilon_4&=0,\quad x\in \Gamma^{wa}_\varepsilon,\;t \in (0,T)\\
-n^\varepsilon\cdot d^\varepsilon_4\nabla u^\varepsilon_4&=0,\quad x\in \Gamma^{sw}_\varepsilon,\;t \in (0,T).
 \end{aligned}
\end{equation}
The mass-balance equation for the gypsum produced at the water-solid interface is
\begin{equation}\label{v}
\begin{aligned}
\partial_t u^\varepsilon_5&=
\eta(u^\varepsilon_1,u^\varepsilon_5),\quad x\in \Gamma^{sw}_\varepsilon, \;t \in (0,T)\\
 u^\varepsilon_5(x,0)&={u^\varepsilon_5}_0(x),\quad x\in \Gamma^{sw}_\varepsilon, \;t \in (0,T).
 \end{aligned}
\end{equation}
\vspace{-0.3in}
\section{Weak formulation and basic results}\label{basic}
\vspace{-0.3in}
We begin this section with a list of notations and function spaces. Then we indicate our working assumptions and give the  weak formulation of the microscopic problem; we bring reader's attention to the well-posedness of the system (\ref{i})--(\ref{v}).
\vspace{-0.3in}
\subsection{Notations and function spaces}
\vspace{-0.3in}
We use  $(\alpha,\beta)_{(0,T) \times
\Omega^\varepsilon}:=\int_0^T\int_{\Omega^\varepsilon}\alpha\beta
 dxdt$,
  $(\alpha,\beta)_{(0,T)\times\Gamma_\varepsilon}:=\int_0^T\int_{\Gamma_\varepsilon}\alpha\beta d\sigma_x
  dt$.
 $\langle\cdot \rangle$, $|\cdot|$ and $\|\cdot\|$ denote the dual pairing of $H^1(\Omega^\varepsilon)$ and $H^{-1}(\Omega^\varepsilon)$,
 the norm in $L^2(\Omega^\varepsilon)$, and the norm in $H^1(\Omega^\varepsilon),$ respectively.
 $\varphi^+$ and $\varphi^-$ will point out the positive and respectively the negative part of the function
 $\varphi$. We denote by $C^\infty_\#(Y)$, $H^1_\#(Y)$, and $H^1_\#(Y)/\mathbb{R}$, the space of infinitely differentiable functions in $\mathbb{R}^n$ that are periodic of period $Y$,
the completion of $C^\infty_\#(Y)$ with respect to $H^1-$norm, and
the respective quotient space,
respectively. Furthermore, $H^1_{\Gamma^D}(\Omega):=\{u\in
H^1(\Omega)|u=0 \mbox{ on } \Gamma^D\}$. The Sobolev space
$H^\beta(\Omega)$ as a completion of
$C^\infty(\Omega)$ is a Hilbert space equipped with a norm
\vspace{-0.2in}
$$\| \varphi\|_{H^\beta(\Omega)}= \| \varphi\|_{H^{[\beta]}(\Omega)}+ \left(\int_\Omega\int_\Omega
\frac{|\varphi(x)-\varphi(y)|^2}{|x-y|^{n+2(\beta-[\beta])}}dxdy\right)^\frac{1}{2}
$$
 and (cf. Theorem 7.57 in \cite{kufner}) the embedding $H^\beta(\Omega)\hookrightarrow L^2(\Omega)$ is continuous.
Since we deal with an evolution problem, we need typical Bochner
spaces like $L^2(0,T;H^1(\Omega))$, $L^2(0,T;H^{-1}(\Omega))$,
$L^2(0,T;H^1_{\Gamma^D}(\Omega))$, and
$L^2((0,T)\times\Omega;H^1_\#(Y)/\mathbb{R})$.
 In the analysis of
the microscopic model, we use frequently the following trace
inequality for $\varepsilon-$dependent hypersurfaces
$\Gamma^{wa}_\varepsilon$: For $\varphi_\varepsilon \in
H^1(\Omega^\varepsilon)$, there exists a constant $C^*$, which is
independent of $\varepsilon$, such that
\begin{equation}\label{1a}
\varepsilon|\varphi_\varepsilon|^2_{L^2(\Gamma_\varepsilon)}\leq C^*
(|\varphi_\varepsilon|^2_{L^2(\Omega^\varepsilon)}+\varepsilon^2 |
\nabla \varphi_\varepsilon|^2_{L^2(\Omega^\varepsilon)}).
\end{equation}
The proof of (\ref{1a}) is given in Lemma 3 of \cite{hornung}. For a
function $\varphi^\varepsilon \in H^\beta(\Omega^\varepsilon)$ with
$\beta\in(\frac{1}{2},1)$, the inequality (\ref{1a}) refines into
\vspace{-0.2in}
\begin{equation}\label{1b}
\varepsilon|\varphi_\varepsilon|^2_{L^2(\Gamma_\varepsilon)}\leq
C^*_0
(|\varphi_\varepsilon|^2_{L^2(\Omega^\varepsilon)}+\varepsilon^{2\beta}
\int_\Omega^\varepsilon\int_\Omega^\varepsilon
\frac{|\varphi^\varepsilon(x)-\varphi^\varepsilon(y)|^2}{|x-y|^{n+2\beta}}dxdy),
\end{equation}
where $C^*_0$ is again a constant independent of $\varepsilon$. For
proof of (\ref{1b}), see \cite{MM}.
To simplify the writing of some of the estimates, we employ the next
set of notations:
 \begin{equation*}\label{}
\begin{aligned}
 d_i&:=\min_{[0,T]\times \bar{\Omega}}\mid d_i^\varepsilon\mid, \;i \in\{1,2,3,4\},\quad\tilde{d}_i:=\min_{[0,T]\times \bar{\Omega}}\mid \tilde{d}_i^\varepsilon\mid, \\
 D_m&:=\max_{[0,T] \times \Omega^\varepsilon}|\partial_t d^\varepsilon_m|, \;m\in\{1,2,3\},\quad
 k_j:=\min_{[0,T]\times \bar{\Omega}} \mid k_j^\varepsilon\mid,\;j\in\{1,2\}\\
K_j&:=\min_{[0,T]\times \bar{\Omega}} \mid \partial_tk_j^\varepsilon\mid,\quad
 \tilde{k}_j:=\min_{[0,T]\times \bar{\Omega}} \mid{\tilde{k}}_j^\varepsilon\mid,\\
k^\infty_m&:= \sup_{(0,T)\times \Omega}\mid  k_{m}^\varepsilon\mid,\quad
\tilde{k}^\infty_m:= \sup_{(0,T)\times \Omega}\mid  \tilde{k}_{m}^\varepsilon\mid,\\
K^\infty_m&:= \sup_{(0,T)\times \Omega}\mid \partial_t k_{m}^\varepsilon\mid,\quad M_i:=\sup_{(0,T)\times \Omega}\mid  u_{i}^\varepsilon\mid,\;i\in\{1,2,3,4,5\},
\end{aligned}
\end{equation*}
\begin{equation*}
\begin{aligned}
A^\infty&:=\sup_{(0,T)\times \Gamma^{wa}_{\varepsilon}}|a_\varepsilon|,\; B^\infty:=\sup_{(0,T)\times \Gamma^{wa}_{\varepsilon}}|b_\varepsilon|,\\
A^\infty&:=\sup_{(0,T)\times \Gamma^{wa}_{\varepsilon}}|\partial_ta_\varepsilon|,\; B^\infty:=\sup_{(0,T)\times \Gamma^{wa}_{\varepsilon}}|\partial_tb_\varepsilon|,\\
\tilde{a}^\infty&:=\sup_{(0,T)\times \Gamma^{wa}}|\tilde{a}|,\; \tilde{b}^\infty:=\sup_{(0,T)\times \Gamma^{wa}}|\tilde{b}|,\\
Q^\infty&:=\sup_{s\in(0,T)\times \Gamma^{sw}_{\varepsilon}}|Q(s)|,\quad
\bar{\eta}:=||\eta||_\infty,\quad
\hat{\eta}:=||\partial_t\eta||_\infty.
 \end{aligned}
\end{equation*}
\vspace{-0.3in}
\subsection{Assumptions on the data and parameters}
\vspace{-0.1in}
We consider the following restriction on the data and parameters:
\begin{enumerate}
\item[(A1)] $d_i\in L^\infty((0,T)\times Y)^{3\times 3}$, $\partial_td_i\in L^\infty((0,T)\times Y)^{3\times 3}$,
 $\partial_{tt}d_i\in L^\infty((0,T)\times Y)^{3\times 3}$, $(d_i(t,x)\xi,\xi)\geq d_{i0}\mid\xi\mid^2$
  for $d_{i0}>0$, for every $\xi\in\mathbb{R}^3$, $(t,x)\in (0,T)\times Y$, $i\in \{1,2,3,4\}$.
\item[(A2)]$\eta$ is measurable w.r.t. $t$ and $x$ and
$\eta(\alpha,\beta)=k_3^\varepsilon R(\alpha)Q(\beta)$, $R$ is sub-linear and locally Lipschitz function
   and $Q$ is bounded and locally Lipschitz function such that
   \begin{eqnarray}
R(\alpha)=
\left \{ \begin{array}{ccc}
\mbox{positive},\;\;\mbox{if}\;\;\alpha\geq0,\\\nonumber
0,\;\;\;\mbox{otherwise}\nonumber
 \end{array} \right.
\quad\quad
Q(\beta)=
\left \{ \begin{array}{ccc}
\mbox{positive},\;\;\mbox{if}\;\;\beta<\beta_{max},\\\nonumber
0,\;\;\;\mbox{otherwise}\nonumber
 \end{array} \right.
\end{eqnarray}
Additionally to (A2), we sometimes assume (A2)', that is
\item[(A2)']$\partial_t\eta\leq\hat{\eta}$.
\item[(A3)] $u^\varepsilon_{i0}\in L^2(\Omega^\varepsilon)\cap L_+^\infty(\Omega^\varepsilon),\;i\in\{1,2,4\}$, $u^\varepsilon_{30}\in  L^2(\Omega^\varepsilon_1)\cap L_+^\infty(\Omega^\varepsilon_1)$,
$u^\varepsilon_{50}\in  L^2(\Gamma^{sw}_\varepsilon)\cap L_+^\infty(\Gamma^{sw}_\varepsilon)$.
\item[(A4)] $a^\infty M_3=b^\infty M_2$, $k^\infty_1 M_1= M_4$, $k_1 M_1=k^\infty_2 M_2$.
\item[(A5)] $a,b\in {C}^1([0,T];C^{0,\alpha}(\Gamma^{wa})), a,b\geq0 $ in $[0,T]\times\Gamma^{wa},\partial_ta,\partial_tb \in  L^\infty((0,T)\times \Gamma^{wa})$.
\item[(A6)] $\partial_t u_3^D$, $\partial_{tt} u_3^D$ and $\nabla \partial_t u_3^D$ are bounded.
\item[(A7)] $k_3\in C^1([0,T];C^{0,\alpha}(\Gamma^{sw}))$ and $k_j\in C^1([0,T];C^{0,\alpha}(\bar{Y}))$ for any $j\in\{1,2\}$ and $\alpha\in ]0,1]$.
\end{enumerate}
The assumptions (A1)--(A3), (A5), and (A6) are of technical nature.
The first equality in (A4) points out an infinitely fast
(equilibrium) Henry law, while the last two equalities remotely
resemble a detailed balance in two of the involved chemical
reactions.
\vspace{-0.3in}
\subsection{Weak formulation of the microscopic model}
\vspace{-0.15in}
\begin{defn}\label{def}
Assume (A1) and (A3). We call the vector
$u^\varepsilon=(u^\varepsilon_1,u^\varepsilon_2,u^\varepsilon_3,u^\varepsilon_4,u^\varepsilon_5)$,
 a weak solution to (\ref{i})--(\ref{v})
 if $u^\varepsilon_j\in L^2(0,T;H^1(\Omega^\varepsilon)), \partial_t u^\varepsilon_j
\in L^2(0,T;H^{-1}(\Omega^\varepsilon)), j \in \{1,2,4\}$, $u^\varepsilon_3\in u_3^D+L^2(0,T;H_{\Gamma^D}^1(\Omega_1^\varepsilon)),
\partial_t u^\varepsilon_3 \in u_3^D+L^2(0,T;H^{-1}(\Omega_1^\varepsilon)),
u^\varepsilon_5\in L^\infty((0,T)\times\Gamma^{sw}_\varepsilon),
\partial_tu^\varepsilon_5\in L^\infty((0,T)\times\Gamma^{sw}_\varepsilon)$ such that
the following identities hold
\begin{equation}\label{def1}
\begin{aligned}
 \langle \partial_t u^\varepsilon_1,\varphi_1\rangle_{(0,T) \times \Omega^\varepsilon}&+
 (d_1\nabla u^\varepsilon_1),\nabla\varphi_1)_{(0,T) \times \Omega^\varepsilon}\\
& =-( k_1u^\varepsilon_1,\varphi_1)_{(0,T) \times \Omega^\varepsilon}
 + ( k_2 u^\varepsilon_2,\varphi_1)_{(0,T) \times \Omega^\varepsilon}\\&-\varepsilon
(\eta (u^\varepsilon_1, u^\varepsilon_4),\varphi_1)_{(0,T) \times \Gamma^{sw}_\varepsilon},
\end{aligned}
\end{equation}
\begin{equation}\label{def2}
\begin{aligned}
  \langle\partial_t u^\varepsilon_2,\varphi_2 \rangle_{(0,T) \times \Omega^\varepsilon}&+(d^\varepsilon_2\nabla u^\varepsilon_2),\nabla\varphi_2)_{(0,T) \times \Omega^\varepsilon}
  \\&=( k^\varepsilon_1u^\varepsilon_1,\varphi_2)_{(0,T) \times \Omega^\varepsilon}
 -( k^\varepsilon_2 u^\varepsilon_2,\varphi_2)_{(0,T) \times \Omega^\varepsilon}\\&+\varepsilon
  (a_\varepsilon u^\varepsilon_3,\varphi_2)_{(0,T) \times \Gamma^{wa}_\varepsilon}-\varepsilon
  (a_\varepsilon u^\varepsilon_2,\varphi_2)_{(0,T) \times \Gamma^{wa}_\varepsilon},
\end{aligned}
\end{equation}
\begin{equation}\label{def3}
\begin{aligned}
  \langle\partial_t u^\varepsilon_3,\varphi_3 \rangle_{(0,T) \times \Omega_1^\varepsilon}&=-
  (d^\varepsilon_3\nabla u^\varepsilon_3),\nabla\varphi_3)_{(0,T) \times \Omega_1^\varepsilon}\\
  &-\varepsilon
  (a_\varepsilon u^\varepsilon_3,\varphi_3)_{(0,T) \times \Gamma^{wa}_\varepsilon}+\varepsilon
  (a_\varepsilon u^\varepsilon_2,\varphi_3)_{(0,T) \times \Gamma^{wa}_\varepsilon},
\end{aligned}
\end{equation}
\begin{equation}\label{def4}
\begin{aligned}
  \langle\partial_t u^\varepsilon_4,\varphi_4 \rangle_{(0,T) \times \Omega^\varepsilon}&=-
  (d^\varepsilon_4\nabla u^\varepsilon_4),\nabla\varphi_4)_{(0,T) \times \Omega^\varepsilon}+
  ( k^\varepsilon_1u^\varepsilon_1,\varphi_4)_{(0,T) \times \Omega^\varepsilon}
\end{aligned}
\end{equation}
for all $\varphi_j \in L^2(0,T;H ^1(\Omega^\varepsilon)), j\in\{1,2,4\}$ and $\varphi_3\in L^2(0,T;H_{\Gamma^D}^1(\Omega_1^\varepsilon))$
 together with the ode
\begin{equation}\label{def5}
\begin{aligned}
\partial_t u^\varepsilon_5&=&
\eta(u^\varepsilon_1,u^\varepsilon_5)\;\;\mbox{a.e. on}\;\;(0,T)\times{\Gamma^{ws}_\varepsilon}
\end{aligned}
\end{equation}
and the initial conditions
\begin{equation}\label{def6}
\begin{aligned}
u^\varepsilon_i(0,x)&=u^\varepsilon_{i0}(x)\;\;x\in \Omega^\varepsilon \mbox{ for all } i\in\{1,2,4\},\\
u^\varepsilon_3(0,x)&=u^\varepsilon_{30}(x)\;\;x\in \Omega_1^\varepsilon,\\
u^\varepsilon_5(0,x)&=u^\varepsilon_{50}(x)\;\;x\in \Gamma^{ws}_\varepsilon.
\end{aligned}
\end{equation}
\end{defn}
\vspace{-0.3in}
\subsection{Basic results}
\vspace{-0.2in}
\begin{lem}\label{Basic results}(Positivity and $L^\infty$-estimates)
Assume (A1)-(A6), and let $t\in [0,T]$ be arbitrarily chosen. Then the following estimates hold:
\vspace{-0.2in}
\begin{itemize}
\item[(i)]$u^\varepsilon_i(t)\geq0,\;\;i\in\{1,2,4\}$ a.e. in ${\Omega^\varepsilon}$, $u^\varepsilon_3(t)\geq0$ a.e. ${\Omega^\varepsilon_1}$ and
 $u^\varepsilon_5(t)\geq0$ a.e. on $\Gamma^{ws}_\varepsilon$.
\item[(ii)] $u^\varepsilon_i(t)\leq M_i$, $i\in\{1,2\}$, $u^\varepsilon_4(t)\leq (t+1)M_4$ a.e. in ${\Omega^\varepsilon}$ , $u^\varepsilon_3(t)\leq M_3$ a.e. in ${\Omega^\varepsilon_1}$ and
$u^\varepsilon_5(t)\leq M_5$ a.e. on $\Gamma^{ws}_\varepsilon$.
\end{itemize}
\end{lem}
\vspace{-0.2in}
\textit{Proof} (i). We test  (\ref{def1})-(\ref{def4}) with
$\varphi=(-{u^\varepsilon_1}^-,-{u^\varepsilon_2}^-,-{u^\varepsilon_3}^-,-{u^\varepsilon_4}^-)$
element of the space $[L^2(0,T;H^1(\Omega^\varepsilon))]^2\times
L^2(0,T; H^1_{\Gamma^D}(\Omega_1^\varepsilon)\times
L^2(0,T;H^1(\Omega^\varepsilon)$. We obtain the following inequality
\vspace{-0.2in}
\begin{equation}\label{z}
\begin{aligned}
\frac{1}{2}\partial_t |{u^\varepsilon_1}^-|^2+d_1|\nabla{u^\varepsilon_1}^-|^2&\leq -k_1|{u^\varepsilon_1}^-|^2+k^\infty_2({u^\varepsilon_1}^-,
{u^\varepsilon_2}^-)\\
&-\varepsilon(\eta (u^\varepsilon_1, u^\varepsilon_5),-{u^\varepsilon_1}^-)_{\Gamma^{sw}_\varepsilon}.
\end{aligned}
\end{equation}
Note that the first term on the r.h.s of (\ref{z}) is negative,
while the third term is zero because of (A2). We then get
\begin{equation}\label{a1}
\partial_t |{u^\varepsilon_1}^-|^2+2d_1|\nabla{u^\varepsilon_1}^-|^2\leq k^\infty_2\left(|{u^\varepsilon_1}^-|^2+|{u^\varepsilon_2}^-|^2\right).
\end{equation}
On the other hand, (\ref{def2}) leads to
\begin{equation*}
\begin{aligned}
\frac{1}{2}\partial_t |{u^\varepsilon_2}^-|^2+d_2|\nabla{u^\varepsilon_2}^-|^2 &\leq \frac{k^\infty_1}{2}\left(|{u^\varepsilon_1}^-|^2
+|{u^\varepsilon_2}^-|^2\right)\nonumber
\\&+\varepsilon a^\infty({u^\varepsilon_2}^-,{u^\varepsilon_3}^-)_{\Gamma^{wa}_\varepsilon}+
\varepsilon b^\infty|{u^\varepsilon_2}^-|^2_{\Gamma^{wa}_\varepsilon}.
\end{aligned}
\end{equation*}
By the trace inequality (\ref{1a}) (with $\varepsilon<1$), we get
\begin{equation}
\begin{aligned}
\partial_t |{u^\varepsilon_2}^-|^2&+2(d_2-C^*b^\infty)|\nabla{u^\varepsilon_2}^-|^2 \leq {k^\infty_1}
\left(|{u^\varepsilon_1}^-|^2+|{u^\varepsilon_2}^-|^2\right)
\\&+2C^*b^\infty|{u^\varepsilon_2}^-|^2+2\varepsilon a^\infty({u^\varepsilon_2}^-,{u^\varepsilon_3}^-)_{\Gamma^{wa}_\varepsilon}.
\end{aligned}
\end{equation}
(\ref{def3}) leads to
\begin{equation}
\partial_t |{u^\varepsilon_3}^-|^2+2(d_3-C^*a^\infty)|\nabla{u^\varepsilon_3}^-|^2\leq 2\varepsilon b^\infty({u^\varepsilon_2}^-,{u^\varepsilon_3}^-)_{\Gamma^{wa}_\varepsilon}+2C^*a^\infty|{u^\varepsilon_3}^-|^2,
\end{equation}
while from (\ref{def4}), we see that
\begin{equation}\label{a2}
\partial_t |{u^\varepsilon_4}^-|^2+2d_4|\nabla{u^\varepsilon_5}^-|^2\leq
k^\infty_1\left(|{u^\varepsilon_1}^-|^2+|{u^\varepsilon_5}^-|^2\right).
\end{equation}
Adding up inequalities (\ref{a1})-(\ref{a2}) gives
\begin{equation}\label{}
\begin{aligned}
\partial_t\sum_{i=1}^4 |{u^\varepsilon_i}^-|^2&+2d_1|\nabla{u^\varepsilon_1}^-|^2+2(d_2-C^*b^\infty)|\nabla{u^\varepsilon_2}^-|^2\\&+
2(d_3-C^*a^\infty)|\nabla{u^\varepsilon_3}^-|^2+
2d_4|\nabla{u^\varepsilon_4}^-|^2\\&\leq
\left(2k^\infty_1+k^\infty_2+2C^*b^\infty+2C^*a^\infty\right)\sum_{{i=1}}^4 |{u^\varepsilon_i}^-|^2\\&+2\varepsilon(a^\infty+b^\infty)
({u^\varepsilon_2}^-,{u^\varepsilon_3}^-)_{\Gamma^{wa}_\varepsilon},
\end{aligned}
\end{equation}
and hence,
\begin{equation}\label{a3}
\begin{aligned}
\partial_t\sum_{i=1}^4 |{u^\varepsilon_i}^-|^2&+2d_1|\nabla{u^\varepsilon_1}^-|^2+2(d_2-C^*b^\infty)|\nabla{u^\varepsilon_2}^-|^2\\&+
2(d_3-C^*a^\infty)|\nabla{u^\varepsilon_3}^-|^2+
2d_4|\nabla{u^\varepsilon_5}^-|^2\\&\leq
\left(2k^\infty_1+k^\infty_2+C^*(a^\infty+b^\infty)\right)\sum_{i=1}^4 |{u^\varepsilon_i}^-|^2\\&+
\varepsilon\left(a^\infty+b^\infty)(\delta|{u^\varepsilon_2}^-|_{\Gamma^{wa}_\varepsilon}^2+
\frac{1}{\delta}|{u^\varepsilon_3}^-|^2_{\Gamma^{wa}_\varepsilon}\right).
\end{aligned}
\end{equation}
Applying the trace inequality (\ref{1a}) to estimate the last term
on the right side of (\ref{a3}), we finally get
\begin{equation*}
\begin{aligned}
\partial_t\sum_{i=1}^4 |{u^\varepsilon_i}^-|^2&+2d_1|\nabla{u^\varepsilon_1}^-|^2+(2d_2-2C^*b^\infty-C^*\delta(a^\infty+b^\infty))
|\nabla{u^\varepsilon_2}^-|^2\\&+
(2d_3-2C^*a^\infty-\frac{C^*2}{\delta}(a^\infty+b^\infty))|\nabla{u^\varepsilon_3}^-|^2+
2d_4|\nabla{u^\varepsilon_4}^-|^2\\&\leq
C_1\sum_{i=1}^4 |{u^\varepsilon_i}^-|^2.
\end{aligned}
\end{equation*}
Thus, we have
\begin{equation*}
\begin{aligned}
\partial_t\sum_{i=1}^4 |{u^\varepsilon_i}^-|^2&\leq&
C_1\sum_{i=1}^4 |{u^\varepsilon_i}^-|^2.
\end{aligned}
\end{equation*}
where
$C_1:=2k^\infty_1+k^\infty_2+C^*(a^\infty+b^\infty)+C^*(\delta+\frac{1}{\delta})(a^\infty+b^\infty)$
and $\delta$ is chosen conveniently. Gronwall's inequality together
with $[u^\varepsilon_i(0)]^-=0$ gives now the desired result. Note
that (A2) ensures automatically the positivity of $u^\varepsilon_5$.

(ii). We consider the test function
\begin{equation*}
(\varphi_1,\varphi_2,\varphi_3,\varphi_4)=((u^\varepsilon_1-M_1)^+,(u^\varepsilon_2-M_2)^+,(u^\varepsilon_3-M_3)^+,(u^\varepsilon_4-(t+1)M_4)^+).
\end{equation*}
Obviously,
$\varphi\in [L^2(0,T;H^1(\Omega^\varepsilon))]^2\times L^2(0,T;H^1_{\Gamma^D}(\Omega_1^\varepsilon)\times
L^2(0,T;H^1(\Omega^\varepsilon)$ is allowed as test function. We obtain from (\ref{def1}) that
\vspace{-0.2in}
\begin{equation*}
\begin{aligned}
\frac{1}{2}\partial_t |(u^\varepsilon_1-M_1)^+|^2&+d_1|\nabla(u^\varepsilon_1-M_1)^+|^2\leq -k_1|(u^\varepsilon_1-M_1)^+|^2\\
&-(k_1M_1, (u^\varepsilon_1-M_1)^+)\\\nonumber&+k^\infty_2((u^\varepsilon_1-M_1)^+,
(u^\varepsilon_2-M_2)^+)\\\nonumber&+(k^\infty_2M_2, (u^\varepsilon_1-M_1)^+)\\\nonumber
&-\varepsilon(\eta (u^\varepsilon_1, u^\varepsilon_5),(u^\varepsilon_1-M_1)^+)_{\Gamma^{sw}_\varepsilon}.
\end{aligned}
\end{equation*}
Relying on (A4), we get the estimate
\vspace{-0.2in}
\begin{eqnarray}\label{a4}
\begin{array}{ccc}
\partial_t |(u^\varepsilon_1-M_1)^+|^2&\leq &k^\infty_2(|(u^\varepsilon_1-M_1)^+|^2+|
(u^\varepsilon_2-M_2)^+|^2).
\end{array}
\end{eqnarray}
(\ref{def2}) in combination with (A4) gives that
\begin{equation}\label{a5}
\begin{aligned}
\partial_t |(u^\varepsilon_2-M_2)^+|^2&+2(d_2-C^*b^\infty)|\nabla(u^\varepsilon_2-M_2)^+|^2\\ &\leq
k^\infty_1(|(u^\varepsilon_1-M_1)^+|^2+|(u^\varepsilon_2-M_2)^+|^2)\\&+2C^*b^\infty|(u^\varepsilon_2-M_2)^+|^2\\
&+2\varepsilon a^\infty ((u^\varepsilon_2-M_2)^+,(u^\varepsilon_3-M_3)^+)_{\Gamma^{wa}_\varepsilon}.
\end{aligned}
\end{equation}
By (\ref{def3}), we obtain
\begin{equation}\label{a6}
\begin{aligned}
\partial_t |(u^\varepsilon_3-M_3)^+|^2&+2(d_3-C^*a^\infty)|\nabla(u^\varepsilon_3-M_3)^+|^2\\ &\leq
2C^*a^\infty|\nabla(u^\varepsilon_3-M_3)^+|^2\\&+
2\varepsilon b^\infty ((u^\varepsilon_2-M_2)^+,(u^\varepsilon_3-M_3)^+)_{\Gamma^{wa}_\varepsilon}.
\end{aligned}
\end{equation}
Using again (A4), (\ref{def4}) yields
\vspace{-0.2in}
\begin{eqnarray}\label{a7}
\begin{array}{ccc}
\partial_t |(u^\varepsilon_4-(t+1)M_4)^+|^2&\leq &k^\infty_1(|(u^\varepsilon_1-M_1)^+|^2+|
(u^\varepsilon_4-(t+1)M_4)^+|^2).
\end{array}
\end{eqnarray}
 Adding up (\ref{a4})--(\ref{a7}) side by side, we get
\begin{equation}\label{221}
\begin{aligned}
\sum_{j=1}^3\partial_t |(u^\varepsilon_j-M_j)^+|^2&+\partial_t |(u^\varepsilon_4-(t+1)M_4)^+|^2+
(2d_2-2C^*b^\infty)|\nabla(u^\varepsilon_2-M_2)^+|^2\\\nonumber
&+(2d_3-2C^*a^\infty)|\nabla(u^\varepsilon_3-M_3)^+|^2 \nonumber\\&\leq
(2k^\infty_2+k^\infty_1+2C^*a^\infty+2C^*b^\infty)(\sum_{j=1}^3 |(u^\varepsilon_j-M_j)^+|^2\nonumber\\&+|(u^\varepsilon_4-(t+1)M_4)^+|^2)
+\varepsilon (a^\infty+b^\infty) (\delta|(u^\varepsilon_2-M_2)^+|^2_{\Gamma^{wa}_\varepsilon}\nonumber\\&+\frac{1}{\delta}|(u^\varepsilon_3-M_3)^+|^2_{\Gamma^{wa}_\varepsilon}).
\end{aligned}
\end{equation}
We use the trace inequality (\ref{1a}) (with $\varepsilon<1$) to deal with the boundary terms in (\ref{221}). Then Gronwall's inequality yields for all $t \in(0,T)$ the following estimate
\vspace{-0.2in}
\begin{equation*}\label{}
\begin{aligned}
u^\varepsilon_j(t)&\leq M_j,\;\;j\in\{1,2,5\}\;a.\;e. \;in \;\Omega^\varepsilon\\
u^\varepsilon_3(t)&\leq M_3,\;a.\;e. \;in \;\Omega_1^\varepsilon\\
u^\varepsilon_4&\leq (t+1)M_4\; \mbox{a.e. in}\;\Omega^\varepsilon.
\end{aligned}
\end{equation*}
Furthermore, by (A2) $u^\varepsilon_5$ is bounded.
\begin{prop}(Uniqueness)\label{uni}
Assume (A1)--(A4). Then there exists at most one weak solution in
the sense of Definition \ref{def}. 
\end{prop}
{\em Proof.} We assume that
${u}^{j,\varepsilon}=(u_1^{j,\varepsilon},u_2^{j,\varepsilon},u_3^{j,\varepsilon},u_4^{j,\varepsilon},u_5^{j,\varepsilon}),
 j\in\{1,2\}$  are two distinct weak solutions in the sense of Definition \ref{def}.
 We set $u_i^{\varepsilon}:=u_i^{1,\varepsilon}-u_i^{2,\varepsilon}$ for all $i\in\{1,2,3,4\}$. Firstly, we deal with (15). We obtain
 \begin{equation}\label{odeinteg}
\begin{aligned}
\partial_t {u^{1,\varepsilon}_5}-\partial_t{u^{2,\varepsilon}_5}=\eta(u^{1,\varepsilon}_1,u^{1,\varepsilon}_5
)-\eta(u^{2,\varepsilon}_1,u^{2,\varepsilon}_5).
\end{aligned}
\end{equation}
Integrating (\ref{odeinteg}) along (0,T) and using (A2), we get
 \begin{equation*}
\begin{aligned}
|u^{1,\varepsilon}_5-u^{2,\varepsilon}_5| \leq  k_3^\infty c_Rc_Q M_1
\int
_0^{t}|{u^{1,\varepsilon}_5}-{u^{2,\varepsilon}_5}| d\tau+k_3^\infty c_RQ^{\infty}\int_0^{t}|{u^{1,\varepsilon}_1}-{u^{2,\varepsilon}_1}| d\tau.
\end{aligned}
\end{equation*}
Gronwall's inequality implies
 \begin{equation}\label{u5}
\begin{aligned}
|{u^{1,\varepsilon}_5}(t)-{u^{2,\varepsilon}_5(t)}|\leq C_2\int_0^t |{u^{1,\varepsilon}_1}-{u^{2,\varepsilon}_1}|d\tau \;\;\mbox{for a.e. }t\in(0,T),
\end{aligned}
\end{equation}
where $C_2:=k_3^\infty c_RQ^{\infty}(1+C_3te^{C_3t})$ and $C_3:=k_3^\infty c_Rc_QM_1$.
We calculate
\begin{equation}\label{222}
\frac{1}{2}\partial_t |{u^\varepsilon_1}|^2+d_1|\nabla{u^\varepsilon_1}|^2\leq -k_1|{u^\varepsilon_1}|^2+k^\infty_2({u^\varepsilon_1},{u^\varepsilon_2})
+\varepsilon(\eta _1-\eta _2,{u^{\varepsilon}_1})_{\Gamma^{sw}_\varepsilon},
\end{equation}
where we denote $\eta _1-\eta _2:=\eta(u^{1,\varepsilon}_1,u^{1,\varepsilon}_5
)-\eta(u^{2,\varepsilon}_1,u^{2,\varepsilon}_5)$. We can write
\begin{equation}\label{5}
\begin{aligned}
\frac{1}{2}\partial_t |{u^\varepsilon_1}|^2+d_1|\nabla{u^\varepsilon_1}|^2&\leq -k_1|{u^\varepsilon_1}|^2+\frac{k^\infty_2}{2}(|{u^\varepsilon_1}|^2+|{u^\varepsilon_2}|^2)
\\&+\varepsilon C_3
({u^{1,\varepsilon}_5}-{u^{2,\varepsilon}_5},{u^{\varepsilon}_1})_{\Gamma^{sw}_\varepsilon}\\&+
\varepsilon k_3^\infty c_R Q^\infty({u^{1,\varepsilon}_1}-{u^{2,\varepsilon}_1},{u^{\varepsilon}_1})_{\Gamma^{sw}_\varepsilon}.
\end{aligned}
\end{equation}
Now, inserting (\ref{u5}) in (\ref{5}) yields
\begin{equation}\label{4}
\begin{aligned}
\frac{1}{2}\partial_t |{u^\varepsilon_1}|^2+d_1|\nabla{u^\varepsilon_1}|^2&\leq-k_1|{u^\varepsilon_1}|^2+\frac{k^\infty_2}{2}(|{u^\varepsilon_1}|^2
+|{u^\varepsilon_2}|^2)
\\&+C_4\varepsilon |{u^\varepsilon_1}|^2_{\Gamma^{sw}_\varepsilon}+\frac{\varepsilon C_3^2}{2\delta}\int_0^t |{u^{\varepsilon}_1}|^2_{\Gamma^{sw}_\varepsilon}d\tau,
\end{aligned}
\end{equation}
where $C_4:= k_3^\infty c_R Q^\infty+\frac{ C_3}{2\delta}$. Using
(\ref{1a}), we estimate the last two terms in (\ref{4}) to obtain
the inequality
\begin{equation}\label{11a}
\begin{aligned}
\frac{1}{2}\partial_t |{u^\varepsilon_1}|^2+d_1|\nabla{u^\varepsilon_1}|^2
&\leq-k_1|{u^\varepsilon_1}|^2+\frac{k^\infty_2}{2}(|{u^\varepsilon_1}|^2+|{u^\varepsilon_2}|^2)+
C^*C_4(|{u^\varepsilon_1}|^2\\&+\varepsilon^2|\nabla{u^\varepsilon_1}|)+C^*\frac{ C_3^2}{2\delta}\int_0^t (|{u^{\varepsilon}_1}|^2+\varepsilon^2|\nabla{u^\varepsilon_1}|^2)d\tau.
\end{aligned}
\end{equation}
Note that  the constant $C^*$, arising from in (\ref{11a}), stems
from (\ref{1a}). Rearranging now the terms, we have
\begin{equation}\label{1}
\begin{aligned}
\partial_t |{u^\varepsilon_1}|^2&+(2d_1-2C^*C_4\varepsilon^2)|\nabla{u^\varepsilon_1}|^2
+2k_1|{u^\varepsilon_1}|^2\leq ({k^\infty_2}+C^*C_4)(|{u^\varepsilon_1}|^2\\&+|{u^\varepsilon_2}|^2)+
C^*\frac{ C_3^2}{2\delta}\int_0^t (|{u^{\varepsilon}_1}|^2+\varepsilon^2|\nabla{u^\varepsilon_1}|^2)d\tau.
\end{aligned}
\end{equation}
Following the same line of arguments as before, we obtain from
(\ref{def2}) that
\begin{equation}
\begin{aligned}
\partial_t |{u^\varepsilon_2}|^2+2d_2|\nabla{u^\varepsilon_2}|^2&\leq
-2k_2|{u^\varepsilon_2}|^2+{k^\infty_1}(|{u^\varepsilon_1}|^2+{u^\varepsilon_2}|^2)\\&+2\varepsilon a^\infty({u^\varepsilon_3},{u^\varepsilon_2})_{\Gamma^{wa}_\varepsilon}+
2\varepsilon b^\infty |{u^{\varepsilon}_2}|^2_{\Gamma^{wa}_\varepsilon},
\end{aligned}
\end{equation}
\vspace{-0.1in}
while from (\ref{def3}), we deduce
\begin{equation}
\begin{aligned}
\partial_t |{u^\varepsilon_3}|^2+2d_3|\nabla{u^\varepsilon_3}|^2&\leq
2\varepsilon b^\infty({u^\varepsilon_2},{u^\varepsilon_3})_{\Gamma^{wa}_\varepsilon}+
2\varepsilon a^\infty |{u^{\varepsilon}_3}|^2_{\Gamma^{wa}_\varepsilon}.
\end{aligned}
\end{equation}
Proceeding similarly, (\ref{def4}) yields
\begin{equation}\label{2}
\begin{aligned}
\partial_t |{u^\varepsilon_4}|^2+2d_4|\nabla{u^\varepsilon_4}|^2&\leq
{k^\infty_2}(|{u^\varepsilon_1}|^2+|{u^\varepsilon_4}|^2).
\end{aligned}
\end{equation}
Putting together (\ref{1})--(\ref{2}), we get
\begin{equation}\label{12a}
\begin{aligned}
\partial_t \Sigma^4_{i=1}|{u^\varepsilon_i}|^2&+(2d_1-C^*C_4\varepsilon^2)|\nabla{u^\varepsilon_1}|^2+2d_2|\nabla{u^\varepsilon_2}|^2
+2d_3|\nabla{u^\varepsilon_3}|^2\\&+2d_4|\nabla{u^\varepsilon_4}|^2+2k_1|{u^\varepsilon_1}|^2\\&\leq
 \left(2{k^\infty_1}+{k^\infty_2}+C^*C_2\right)\Sigma^4_{i=1}|{u^\varepsilon_i}|^2\\&+C^*\frac{ C_1^2}{2\delta}\int_0^t (|{u^{\varepsilon}_1}|^2+\varepsilon^2|\nabla{u^\varepsilon_1}|^2)d\tau\\&+
2\varepsilon b |{u^{\varepsilon}_2}|^2_{\Gamma^{wa}_\varepsilon}+
2\varepsilon a |{u^{\varepsilon}_3}|^2_{\Gamma^{wa}_\varepsilon}\\&+\varepsilon (a^\infty+b^\infty)(\delta|{u^\varepsilon_2}|_{\Gamma^{wa}_\varepsilon}^2+\frac{1}{\delta}|{u^\varepsilon_3}|_{\Gamma^{wa}_\varepsilon}^2).
\end{aligned}
\end{equation}
Applying the trace inequality (\ref{1a}) to the boundary terms in (\ref{12a}), we get
\begin{equation}
\begin{aligned}
\partial_t \Sigma^4_{i=1}|{u^\varepsilon_i}|^2&+(2d_1-2C^*C_4\varepsilon^2)|\nabla{u^\varepsilon_1}|^2\\&+(2d_2-
2C^*b^\infty\varepsilon^2-C^*\delta\varepsilon^2(a^\infty+b^\infty))
|\nabla{u^\varepsilon_2}|^2\\&
+(2d_3-2C^*a^\infty\varepsilon^2-\frac{C^*\varepsilon^2}{\delta}(a^\infty+b^\infty))|\nabla{u^\varepsilon_3}|^2\\
&+2d_4|\nabla{u^\varepsilon_4}|^2+
2k_1|{u^\varepsilon_1}|^2 \leq
 C_5\Sigma^4_{i=1}|{u^\varepsilon_i}|^2\\
 \end{aligned}\nonumber
\end{equation}
\begin{equation}\label{6}
\begin{aligned}
&+C^*\frac{ C_1^2}{2\delta}\int_0^t (|{u^{\varepsilon}_1}|^2+\varepsilon^2|\nabla{u^\varepsilon_1}|^2)d\tau,
\end{aligned}
\end{equation}
where
$C_5:=2{k^\infty_1}+{k^\infty_2}+C^*C_2+2C^*(a^\infty+b^\infty)+C^*(a^\infty+b^\infty)(\delta+\frac{1}{\delta})$.
Let us choose $\varepsilon$ and $\delta$ such that
\begin{equation*}
\begin{aligned}
\varepsilon&\in\left]0,\frac{2d_1}{C_1C^*}\right[\\
\delta &\in\left[\frac{C^*\varepsilon^2(a^\infty+b^\infty)}{2d_3-C^*a^\infty\varepsilon^2},\frac{2d_2-C^*b^\infty\varepsilon^2}{C^*\varepsilon^2(a^\infty+b^\infty)}\right].
\end{aligned}
\end{equation*}
With this choice of $(\varepsilon,\delta)$, (\ref{6}) takes the form
\begin{equation*}
\partial_t \Sigma^4_{i=1}|{u^\varepsilon_i}|^2+\bar{C}|\nabla{u^\varepsilon_1}|^2+
\bar{C}|{u^\varepsilon_1}|^2\leq
 {C_6}{}(\Sigma^4_{i=1}|{u^\varepsilon_i}|^2+\int_0^t (|{u^{\varepsilon}_1}|^2+\varepsilon^2|\nabla{u^\varepsilon_1}|^2)d\tau),
\end{equation*}
where
$C_6:=2{k^\infty_1}+{k^\infty_2}+C^*C_2+C^*(a^\infty+b^\infty)+C^*\frac{
C_1^2}{2\delta}$ and
$\bar{C}:=\min\{2d_1-2C^*C_2\varepsilon^2,2k_1\}$. Taking in
(\ref{6}) the supremum along $t \in(0,T)$ and applying Gronwall's
inequality, we obtain the following estimate
\begin{equation}\label{3a}
 \Sigma^4_{i=1}|{u^\varepsilon_i}|^2+\bar{C}\int_0^T|\nabla{u^\varepsilon_1}|^2dt+
\bar{C}\int_0^T|{u^\varepsilon_1}|^2dt\leq0.
\end{equation}
Thus, the proof of Proposition \ref{uni} is completed.
\vspace{-0.1in}
\begin{thm}(Global Existence)\label{existence}
Assume $(A1)-(A3)$. Then there exists at least a global-in-time weak
solution in the sense of Definition \ref{def}.
\end{thm}

{\em Proof}. The proof is based on the Galerkin argument. Since the
proof is rather standard, and here we wish to focus on the passage
to the limit $\varepsilon\to 0$, we omit it.
\vspace{-0.3in}
\section{{\em A priori} estimates for microscopic solutions}\label{est}
\vspace{-0.15in}
This section includes the $\varepsilon-$ independent estimates.
\begin{lem}\label{lemma2}
Assume (A1)-(A6). Then the weak solution of the microscopic model
(\ref{def1})-(\ref{def6}) satisfies the following  {\em a priori}
bounds:
\vspace{-0.1in}
 \begin{eqnarray}
&&\parallel u^\varepsilon_{j}\parallel_{L^2(0,T;H^1(\Omega^\varepsilon))}\leq C,\;j\in\{1,2,3,4\}\label{abc}\\
&&\parallel \nabla \partial_tu^\varepsilon_{2}\parallel_{L^2(0,T;L^2(\Omega^\varepsilon))}\leq C,\label{ali} \\
&&\parallel\partial_tu^\varepsilon_{j}\parallel_{L^2(0,T;L^2(\Omega^\varepsilon))}\leq C, \label{*}\\
&&\parallel u^\varepsilon_{3}\parallel_{L^2(0,T;H^1(\Omega^\varepsilon_1))}\leq C,\\
&&\parallel \nabla \partial_tu^\varepsilon_{3}\parallel_{L^2(0,T;L^2(\Omega^\varepsilon_1))}\leq C, \label{bb} \\ &&\parallel\partial_tu^\varepsilon_{3}\parallel_{L^2(0,T;L^2(\Omega^\varepsilon_1))}\leq C,
\end{eqnarray}
\begin{eqnarray}
&&\parallel u^\varepsilon_{5}\parallel_{L^\infty((0,T)\times \Gamma^{sw}_\varepsilon)}\leq C,\\
&&\parallel \partial_tu^\varepsilon_{5}\parallel_{L^2((0,T)\times \Gamma^{sw}_\varepsilon)}\leq C.\label{acd}
\end{eqnarray}
In (\ref{abc})--(\ref{acd}), the generic constant $C$ is independent of $\varepsilon$.
\end{lem}
{\em Proof.} We test (\ref{def1}) with $\varphi_1=u^\varepsilon_{1}$
to get
 \begin{equation}\label{ass}
 \begin{aligned}
\frac{1}{2}\partial_t|u^\varepsilon_{1}|^2+d_1|\nabla u^\varepsilon_{1}|^2&\leq-k_1|u^\varepsilon_{1}|^2
+k_2^\infty(u^\varepsilon_{1},u^\varepsilon_{2})-\varepsilon(\eta,u^\varepsilon_{1})_{\Gamma^{sw}_\varepsilon},\\
&\leq\frac{k_2^\infty}{2}(|u^\varepsilon_{1}|^2+|u^\varepsilon_{2}|^2)+\varepsilon k^\infty_3Q^\infty c_R
(u^\varepsilon_{1},u^\varepsilon_{1})_{\Gamma^{sw}_\varepsilon}.
\end{aligned}
\end{equation}
After applying the trace inequality to the last term on r.h.s of (\ref{ass}), we get
\begin{equation*}
\begin{aligned}
\frac{1}{2}\partial_t|u^\varepsilon_{1}|^2+d_1|\nabla u^\varepsilon_{1}|^2
&\leq\frac{k_2^\infty}{2}(|u^\varepsilon_{1}|^2+|u^\varepsilon_{2}|^2)+C^*
k^\infty_3Q^\infty c_R
(|u^\varepsilon_{1}|^2+\varepsilon^2|\nabla u^\varepsilon_{1}|^2)_{\Gamma^{sw}_\varepsilon}.
\end{aligned}
\end{equation*}
\begin{equation}\label{24}
\begin{aligned}
\frac{1}{2}\partial_t|u^\varepsilon_{1}|^2&+(d_1-\varepsilon^2C^*
 k^\infty_3Q^\infty c_R)|\nabla u^\varepsilon_{1}|^2\leq C_7(|u^\varepsilon_{1}|^2+|u^\varepsilon_{2}|^2),
\end{aligned}
\end{equation}
where $C_7:=\frac{k_2^\infty}{2}+C^*
 k^\infty_3Q^\infty c_R.$
Taking $\varphi_2=u^\varepsilon_{2}$ in (\ref{def2}), we get
\begin{equation*}
\begin{aligned}
\frac{1}{2}\partial_t|u^\varepsilon_{2}|^2+d_2|\nabla u^\varepsilon_{2}|^2
&\leq\frac{k_1^\infty}{2}(|u^\varepsilon_{1}|^2+|u^\varepsilon_{2}|^2)-k_2|u^\varepsilon_{2}|^2\\&+
\varepsilon a^\infty(u^\varepsilon_{3},u^\varepsilon_{2})_{\Gamma^{wa}_\varepsilon}+\varepsilon b^\infty
| u^\varepsilon_{2}|^2_{\Gamma^{wa}_\varepsilon}.
\end{aligned}
\end{equation*}
Application of the trace inequality (\ref{1a}) only to the last term leads to
\begin{equation}
\begin{aligned}
\frac{1}{2}\partial_t|u^\varepsilon_{2}|^2+(d_2-C^*b^\infty\varepsilon^2)|\nabla u^\varepsilon_{2}|^2
&\leq\frac{k_1^\infty}{2}(|u^\varepsilon_{1}|^2+|u^\varepsilon_{2}|^2)+2\\&+
\varepsilon a^\infty(u^\varepsilon_{3},u^\varepsilon_{2})_{\Gamma^{wa}_\varepsilon}.
\end{aligned}
\end{equation}
We choose $\varphi_3=u^\varepsilon_{3}$ as a test function in (\ref{def3}) to calculate
\vspace{-0.2in}
\begin{equation}
\begin{aligned}
\frac{1}{2}\partial_t|u^\varepsilon_{3}|^2+(d_3-C^*a^\infty\varepsilon^2)|\nabla u^\varepsilon_{3}|^2
&\leq\varepsilon b^\infty(u^\varepsilon_{3},u^\varepsilon_{2})_{\Gamma^{wa}_\varepsilon}+C^*a^\infty|u^\varepsilon_{3}|^2.
\end{aligned}
\end{equation}
\vspace{-0.2in}
Setting $\varphi_4=u^\varepsilon_{4}$ in (\ref{def4}), we are led to
\vspace{-0.2in}
\begin{eqnarray}\label{25}
\frac{1}{2}\partial_t|u^\varepsilon_{4}|^2+d_4|\nabla u^\varepsilon_{4}|^2
&\leq&\frac{k_1^\infty}{2}(|u^\varepsilon_{1}|^2+|u^\varepsilon_{4}|^2).
\end{eqnarray}
Putting together (\ref{24})-(\ref{25}), we obtain
\begin{equation}\label{26}
\begin{aligned}
\frac{1}{2}\Sigma_{i=1}^4\partial_t|u^\varepsilon_{i}|^2&+(d_1-\varepsilon^2C^*
 k^\infty_3
 Q^\infty c_R)|\nabla u^\varepsilon_{1}|^2+d_4|\nabla u^\varepsilon_{4}|^2\\&+
 (d_2-C^*b^\infty\varepsilon^2)|\nabla u^\varepsilon_{2}|^2+(d_3-C^*a^\infty\varepsilon^2)|\nabla u^\varepsilon_{3}|^2
\\&\leq({k_1^\infty}{}+\frac{k_2^\infty}{2}+C^*b^\infty+C^*a^\infty)\Sigma_{i=1}^4|u^\varepsilon_{i}|^2\\&+
\varepsilon(a^\infty+ b^\infty)(u^\varepsilon_{3},u^\varepsilon_{2})_{\Gamma^{wa}_\varepsilon}.
\end{aligned}
\end{equation}
Combing Young's inequality and the trace inequality to the boundary term, (\ref{26}) turns out to be
\begin{equation*}
\begin{aligned}
\frac{1}{2}\Sigma_{i=1}^4\partial_t|u^\varepsilon_{i}|^2&+(d_1-\varepsilon^2C^*
k^\infty_3Q^\infty c_R)|\nabla u^\varepsilon_{1}|^2\\ &+
 (d_2-C^*b\varepsilon^2-\frac{C^*\varepsilon^2\delta}{2}(a^\infty+ b^\infty))|\nabla u^\varepsilon_{2}|^2
\\&+(d_3-C^*a\varepsilon^2-\frac{C^*\varepsilon^2}{2\delta}(a^\infty+ b^\infty))|\nabla u^\varepsilon_{3}|^2
+d_4|\nabla u^\varepsilon_{4}|^2
\\&\leq({k_1^\infty}{}+\frac{k_2^\infty}{2}+C^*(a^\infty+ b^\infty)(\delta+\frac{1}{\delta}))\Sigma_{i=1}^4|u^\varepsilon_{i}|^2.
\end{aligned}
\end{equation*}
Choosing $\varepsilon$ small enough and $\delta$ conveniently such that the coefficients of the terms involving $|\nabla u^\varepsilon_{2}|^2$ and
$|\nabla u^\varepsilon_{3}|^2$ stay positive, we are led to
\vspace{-0.2in}
\begin{eqnarray*}\label{}
\Sigma_{i=1}^4\partial_t|u^\varepsilon_{i}|^2+d'_1|\nabla u^\varepsilon_{1}|^2+d'_2|\nabla u^\varepsilon_{2}|^2+d'_3|\nabla u^\varepsilon_{3}|^2
&+&2d_4|\nabla u^\varepsilon_{4}|^2\leq C_7\Sigma_{i=1}^4|u^\varepsilon_{i}|^2,
\end{eqnarray*}
where
\begin{equation*}
\begin{aligned}
d'_1&:=2(d_1-\varepsilon^2C^* k^\infty_3Q^\infty c_R),\\
 d'_2&:=2(d_2-C^*b^\infty\varepsilon^2-\frac{C^*\varepsilon^2\delta}{2}(a^\infty+ b^\infty)), \\
  d'_3&:=2(d_3-C^*a^\infty\varepsilon^2-\frac{C^*b\varepsilon^2}{2\delta}(a^\infty+ b^\infty)),
\end{aligned}
\end{equation*}
while the constant $C$ is given by
 $$C_8:=2{k_1^\infty}{}+\frac{k_2^\infty}{2}+C^*a^\infty+C^*b^\infty+C^*(a^\infty+ b^\infty)(\delta+\frac{1}{\delta}).$$
Summarizing, we have
\vspace{-0.2in}
 \begin{eqnarray}\label{39}
\Sigma_{i=1}^4\partial_t|u^\varepsilon_{i}|^2+d_0\Sigma_{j=1}^3|\nabla u^\varepsilon_{j}|^2
+d_0|\nabla u^\varepsilon_{3}|^2\leq C\Sigma_{i=1}^4|u^\varepsilon_{i}|^2,
\end{eqnarray}
where $d_0:=min\{d'_1,d'_2,d'_3,d'_4\}$.
By Gronwall's inequality, we have
\vspace{-0.2in}
\begin{eqnarray*}\label{}
\Sigma_{i=1}^4|u^\varepsilon_{i}|^2\leq C\Sigma_{i=1}^4|u_{i}(0)|^2,
\end{eqnarray*}
and hence,
\begin{eqnarray}\label{38}
\parallel u^\varepsilon_{j}\parallel_{L^2(0,T;L^2(\Omega^\varepsilon))}\leq C\mbox{ for all }i\in\{1,2,4\}\;\mbox{and}\;
\|u^\varepsilon_{3}\|_{{L^2(0,T;L^2(\Omega^\varepsilon_1))}}\leq C,
\end{eqnarray}
where $C$ depends on initial data and model parameters but is independent of $\varepsilon$.
Integrating (\ref{39}) along $(0,T)$, we get
\begin{equation}\label{27}
\begin{aligned}
\parallel u^\varepsilon_{j}\parallel_{ L^2(0,T;H^1(\Omega^\varepsilon))}&\leq C, \;j\in\{1,2,4\},\\
\parallel u^\varepsilon_{3}\parallel_{L^2(0,T;H^1(\Omega^\varepsilon_1))}&\leq C.
\end{aligned}
\end{equation}
With the help of (A2) together with the boundedness of $u^\varepsilon_1$, we conclude from (\ref{def5}) that
\begin{equation*}
\begin{aligned}
\parallel u^\varepsilon_{5}\parallel_{L^\infty((0,T)\times\Gamma^{sw}_\varepsilon)}\leq C.
\end{aligned}
\end{equation*}
Multiplying (\ref{def5}) by $\partial_tu^\varepsilon_{5}$ and using (A2),
we get
\begin{equation*}
\begin{aligned}
\|\partial_tu^\varepsilon_{5}\|_{L^2((0,T)\times\Gamma^{sw}_\varepsilon)}\leq C.
\end{aligned}
\end{equation*}
Now, we focus on obtaining $\varepsilon-$independent estimates on the time derivative of the concentrations.
Firstly, we choose $\varphi_1=\partial_tu^\varepsilon_{1}$ and get
\begin{equation}\label{***}
\begin{aligned}
\int_0^t\int_{\Omega^\varepsilon}\partial_tu^\varepsilon_{1}\partial_tu^\varepsilon_{1}dxd\tau&+
\int_0^t\int_{\Omega^\varepsilon}d_1^\varepsilon\nabla u^\varepsilon_{1}\nabla \partial_tu^\varepsilon_{1}dxd\tau\\
&=-\int_0^t\int_{\Omega^\varepsilon}k^\varepsilon_1u^\varepsilon_{1}\partial_tu^\varepsilon_{1}dxd\tau
+\int_0^t\int_{\Omega^\varepsilon}k^\varepsilon_2u^\varepsilon_{2}\partial_tu^\varepsilon_{1}dxd\tau\\&-\varepsilon
\int_0^t\int_{\Gamma^{sw}_\varepsilon}\eta\partial_tu^\varepsilon_{1}d\sigma_xd\tau.
\end{aligned}
\end{equation}
Consequently, it holds
\begin{equation*}
\begin{aligned}
\int_0^t\int_{\Omega^\varepsilon}|\partial_tu^\varepsilon_{1}|^2dxd\tau&+\int_0^t\int_{\Omega^\varepsilon}\left(\frac{1}{2}\partial_t(d_1^\varepsilon|\nabla u^\varepsilon_{1}|^2)-(\partial_td_1^\varepsilon)|\nabla u^\varepsilon_{1}|^2\right)dxd\tau\\&\leq-\frac{k_1}{2}\int_0^t\int_{\Omega^\varepsilon}\partial_t|u^\varepsilon_{1}|^2dxd\tau\\&+\frac{k_2^\infty}{2}
\int_0^t\int_{\Omega^\varepsilon}
\left(\frac{1}{\delta}|u^\varepsilon_{2}|^2+\delta |\partial_tu^\varepsilon_{1}|^2\right)dxd\tau\\&-
{\varepsilon}{}\int_0^t\int_{\Gamma^{sw}_\varepsilon}(\partial_t(\eta u^\varepsilon_{1})-(\partial_t\eta) u^\varepsilon_{1})d\sigma_xd\tau,
\end{aligned}
\end{equation*}
\begin{equation}
\begin{aligned}
(1-\frac{k_2^\infty\delta}{2})\int_0^t\int_{\Omega^\varepsilon}|\partial_tu^\varepsilon_{1}|^2dxd\tau
&\leq D_1\int_0^t\int_{\Omega^\varepsilon}|\nabla u^\varepsilon_{1}|^2dxd\tau\\&+\frac{d_1^\infty}{2}\int_{\Omega^\varepsilon}|\nabla {u_{1}}_0|^2
dx+\frac{k_2^\infty}{2\delta}\int_0^t\int_{\Omega^\varepsilon}
|u^\varepsilon_{2}|^2dxd\tau\\
&+\frac{\varepsilon}{2} \int_{\Gamma^{sw}_\varepsilon}
\left(|\eta|^2+|u^\varepsilon_{1}|^2+|\eta(0)|^2+|u^\varepsilon_{1}(0)|^2\right)d\sigma_x\\
&+
\frac{\varepsilon}{2}\int_0^t \int_{\Gamma^{sw}_\varepsilon}
\left(|\partial_t\eta|^2+|u^\varepsilon_{1}|^2\right)d\sigma_xd\tau,
\end{aligned}
\end{equation}
where $\eta(0):=\eta(u^\varepsilon_{1}(0),u^\varepsilon_{5}(0))$.
Applying (\ref{1a}) and recalling (\ref{27}), we have
\begin{eqnarray}\label{28}
\int_0^t\int_{\Omega^\varepsilon}|\partial_tu^\varepsilon_{1}|^2dxd\tau&\leq&C_9,
\end{eqnarray}
where
\begin{equation*}
\begin{aligned}
C_9&:=D_1\int_0^t\int_{\Omega^\varepsilon}|\nabla u^\varepsilon_{1}|^2dxd\tau+\frac{k_1}{2}\int_{\Omega^\varepsilon}| u^\varepsilon_{1}(0)|^2dx+
\frac{d_1^\infty}{2}\int_{\Omega^\varepsilon}|\nabla {u_{1}}_0|^2
dx
\end{aligned}
\end{equation*}
\begin{equation*}
\begin{aligned}
&+\frac{k_2^\infty}{2\delta}\int_0^t\int_{\Omega^\varepsilon}
|u^\varepsilon_{2}|^2dxd\tau+\frac{\varepsilon}{2}  \int_{\Gamma^{sw}_\varepsilon} (|\overline{\eta}|^2+|{\eta}(0)|^2+|\hat{\eta}|^2)\\
&+
\frac{C^*}{2}\int_0^t\int_{\Omega^\varepsilon}(| u^\varepsilon_{1}|^2+\varepsilon^2|\nabla u^\varepsilon_{1}|^2)dxd\tau+
\frac{C^*}{2}\int_{\Omega^\varepsilon}(| u^\varepsilon_{1}|^2+\varepsilon^2|\nabla u^\varepsilon_{1}|^2)dxd\tau,
\end{aligned}
\end{equation*}
and $\delta \in \left]0, \frac{2}{k_2^\infty}\right[$. Testing
(\ref{def2}) with $\varphi_2=\partial_tu^\varepsilon_{2}$ gives
\vspace{-0.2in}
\begin{eqnarray*}
\int_0^t\int_{\Omega^\varepsilon}|\partial_tu^\varepsilon_{2}|^2dxd\tau&+&\int_0^t\int_{\Omega^\varepsilon}
(\frac{1}{2}\partial_t(d_2^\varepsilon|\nabla u^\varepsilon_{2}|^2)-(\partial_td_2^\varepsilon)|\nabla u^\varepsilon_{2}|^2)dxd\tau\\&\leq&-\frac{k_2}{2}\int_0^t\int_{\Omega^\varepsilon}\partial_t|u^\varepsilon_{2}|^2dxd\tau+\frac{k_1^\infty}{2}
\int_0^t\int_{\Omega^\varepsilon}
(\frac{1}{\delta}|u^\varepsilon_{1}|^2\\&+&\delta |\partial_tu^\varepsilon_{2}|^2)dxd\tau+
\frac{\varepsilon a^\infty}{2}\int_0^t\int_{\Gamma^{wa}_\varepsilon}
\left(|u^\varepsilon_{3}|^2+|\partial_tu^\varepsilon_{2}|^2\right)d\sigma_xd\tau\\&+&
\frac{\varepsilon b^\infty}{2}\int_0^t\int_{\Gamma^{wa}_\varepsilon}\partial_t|u^\varepsilon_{2}|^2d\sigma_xd\tau,
\end{eqnarray*}
\vspace{-0.2in}
and hence,
\vspace{-0.2in}
\begin{eqnarray*}
\int_0^t\int_{\Omega^\varepsilon}|\partial_tu^\varepsilon_{2}|^2dx\tau&+&\frac{d_2}
{2}\int_{\Omega^\varepsilon}|\nabla {u_{2}}^\varepsilon|^2dxd\tau\\&\leq& \frac{d_2^\infty}{2}\int_{\Omega^\varepsilon}|\nabla u^\varepsilon_{2}(0)|^2dx+
D_2\int_0^t\int_{\Omega^\varepsilon}|\nabla u^\varepsilon_{2}|^2dxd\tau\\&+& \frac{k_1^\infty}{2}
\int_0^t\int_{\Omega^\varepsilon}
(\frac{1}{\delta}|u^\varepsilon_{1}|^2+\delta |\partial_tu^\varepsilon_{2}|^2)dxd\tau\\&+&
\frac{C^* a^\infty}{2}\int_0^t\int_{\Omega^\varepsilon}
\left(|u^\varepsilon_{3}|^2+\varepsilon^2|\nabla u^\varepsilon_{3}|^2+\varepsilon^2|\nabla\partial_tu^\varepsilon_{2}|^2\right)dxd\tau\\&+&
\frac{\varepsilon b^\infty}{2}\int_{\Gamma^{wa}_\varepsilon}(|u^\varepsilon_{2}|^2- |u^\varepsilon_{2}(0)|^2)d\sigma_x.
\end{eqnarray*}
By (\ref{1a}) and (\ref{27}), we get
\begin{eqnarray*}
\left(1-\frac{C^*a^\infty}{2}-\frac{k_1^\infty \delta}{2}\right)\int_0^t\int_{\Omega^\varepsilon}|\partial_tu^\varepsilon_{2}|^2dxd\tau
&\leq& C_{10}\left(1+
\varepsilon^2\int_0^t\int_{\Omega^\varepsilon}|\nabla\partial_tu^\varepsilon_{2}|^2dxd\tau\right).
\end{eqnarray*}
Consequently, choosing $\delta
\in]0,\frac{2-C^*a^\infty}{k_1^\infty}[$, we are led to
\begin{equation}\label{29}
\int_0^t\int_{\Omega^\varepsilon}|\partial_tu^\varepsilon_{2}|^2dxd\tau
\leq C_{10}(1+
\varepsilon^2\int_0^t\int_{\Omega^\varepsilon}|\nabla\partial_tu^\varepsilon_{2}|^2dxd\tau),
\end{equation}
where
\begin{eqnarray*}
C_{10}&:=&D_2\int_0^t\int_{\Omega^\varepsilon}|\nabla u^\varepsilon_{1}|^2dxd\tau+
\frac{d_2^\infty}{2}\int_{\Omega^\varepsilon}|\nabla {u^{\varepsilon}}_2(0)|^2
dx+\frac{k_1^\infty}{2\delta}\int_0^t\int_{\Omega^\varepsilon}
|u^\varepsilon_{2}|^2dxd\tau\\&+&\frac{C^*b^\infty}{2} \int_{\Omega^\varepsilon}\left(|u^\varepsilon_{2}|^2+
\varepsilon^2|\nabla u^\varepsilon_{2}|^2+|u^\varepsilon_{2}(0)|^2+\varepsilon^2|\nabla u^\varepsilon_{2}(0)|^2\right)dx\\&+&
\frac{C^*a^\infty}{2} \int_0^t\int_{\Omega^\varepsilon}(|u^\varepsilon_{3}|^2+\varepsilon^2
|\nabla u^\varepsilon_{3}|^2)dxd\tau.
\end{eqnarray*}
The initial data  $u^\varepsilon_{30}$ holding in $\Omega^\varepsilon_1$ and the Dirichlet data $u_3^D$ acting on the exterior
 boundary of $\Omega_1^\varepsilon$ are considered here as  restrictions of the respective
 functions
defined on whole of $\overline{\Omega}$. Testing now (\ref{def3})
with $\varphi_3=\partial_t(u^\varepsilon_{3}-u_3^D)$ leads to
\begin{equation*}\label{}
\begin{aligned}
\int_0^t\int_{\Omega^\varepsilon}|\partial_tu^\varepsilon_{3}|^2dxd\tau&+\frac{d_3}{2}\int_{\Omega^\varepsilon}|\nabla {u^{\varepsilon}}_3|^2\\
&\leq \frac{d_3}{2}\int_{\Omega^\varepsilon}|\nabla {u^\varepsilon_{3}}(0)|^2+\frac{1}{2}(|\partial_tu^\varepsilon_{3}|^2+|\partial_tu^D_{3}|^2)
\\&+D_3\int_0^t\int_{\Omega^\varepsilon}|\nabla {u^\varepsilon_{3}}|^2
+ \frac{d_3^\infty}{2}\int_0^t\int_{\Omega^\varepsilon}(|\nabla {u^\varepsilon_{3}}|^2+|\nabla\partial_tu^D_{3}|^2)
\\&+\frac{{\varepsilon a^\infty}}{\delta}\int_0^t\int_{\Gamma^{wa}_\varepsilon}|u^\varepsilon_{3}|^2
+\frac{\varepsilon\delta}{2} (a^\infty+a^\infty)\int_0^t\int_{\Gamma^{wa}_\varepsilon}|\partial_tu^\varepsilon_{3}|^2
\\&+\frac{\varepsilon}{2} (a^\infty+b^\infty)\int_0^t\int_{\Gamma^{wa}_\varepsilon}|\partial_tu^D_{3}|^2
+\frac{\varepsilon b^{\infty}}{\delta}\int_0^t\int_{\Gamma^{wa}_\varepsilon}|u^\varepsilon_{2}|^2.
\end{aligned}
\end{equation*}
Using (\ref{1a}) and (A6), we obtain
\begin{equation}\label{30}
\int_0^t\int_{\Omega^\varepsilon}|\partial_tu^\varepsilon_{3}|^2dxd\tau\leq C_{11}(1+
\varepsilon^2\delta\int_0^t\int_{\Omega^\varepsilon}|\nabla\partial_tu^\varepsilon_{3}|^2dxd\tau),
\end{equation}
where $\delta \in ]0,\frac{2}{C^*(a^\infty+b^\infty)}[$ and
\begin{equation*}\label{}
\begin{aligned}
C_{11}&:=D_3\int_0^t\int_{\Omega^\varepsilon}|\nabla u^\varepsilon_{3}|^2dxd\tau+\frac{d_3}{2}\int_{\Omega^\varepsilon}|\nabla {u^\varepsilon_{3}}(0)|^2
dx+\frac{1}{2\delta}\int_0^t\int_{\Omega^\varepsilon}|\nabla u^D_{3}|^2\\&+\frac{d^\infty_3}{2}\int_0^t\int_{\Omega^\varepsilon}
(|\nabla u^\varepsilon_{3}|^2+|\nabla \partial_tu^D_{3}|^2)
+\frac{C^*a^\infty }{\delta}\int_0^t\int_{\Omega^\varepsilon}
(| u^\varepsilon_{3}|^2+\varepsilon^2|\nabla u^\varepsilon_{3}|^2)
\\&+\frac{C^*b^\infty}{\delta} \int_0^t\int_{\Omega^\varepsilon}(|u^\varepsilon_{2}|^2+\varepsilon^2
|\nabla u^\varepsilon_{2}|^2)dxd\tau\\&+\frac{C^*(a^\infty+b^\infty)}{2} \int_0^t\int_{\Omega^\varepsilon}(|\partial_tu^D_{3}|^2+\varepsilon^2|\nabla\partial_tu^D_{3}|^2)dxd\tau.
\end{aligned}
\end{equation*}
From (\ref{def4}), we get
\begin{eqnarray}
\int_0^t\int_{\Omega^\varepsilon}|\partial_tu^\varepsilon_{4}|^2dxd\tau
&\leq& C_{12}.
\end{eqnarray}
In order to estimate (\ref{29}) and (\ref{30}), we proceed first
with differentiating ({\ref{def2}}) with respect to time and then
testing the result with $\partial_tu^\varepsilon_{2}$. Consequently,
we derive
\begin{eqnarray}\label{aqw}
\frac{1}{2}\int_{\Omega^\varepsilon}|\partial_tu^\varepsilon_{2}|^2dx
&+&{d_2}
\int_0^t\int_{\Omega^\varepsilon}|\nabla \partial_t{u_{2}}^\varepsilon|^2dxd\tau\\\nonumber
\end{eqnarray}
\begin{eqnarray}
&+&
\int_0^t\int_{\Omega^\varepsilon}\left(\frac{1}{2}(\partial_td_{2}|\nabla {u_{2}}^\varepsilon|^2-(\partial_t\partial_td_{2})|\nabla{u_{2}}^\varepsilon|^2\right)
\nonumber\\
&\leq& \frac{k_1^\infty}{2}\int_0^t\int_{\Omega^\varepsilon}(|\partial_tu^\varepsilon_{1}|^2+|\partial_tu^\varepsilon_{2}|^2)+
\frac{K_1^\infty}{2}\int_0^t\int_{\Omega^\varepsilon}(|u^\varepsilon_{1}|^2+|\partial_tu^\varepsilon_{2}|^2)\nonumber\\
&-&k_2\int_0^t\int_{\Omega^\varepsilon}|\partial_tu^\varepsilon_{2}|^2
-\frac{K_2}{2}\int_0^t\int_{\Omega^\varepsilon}\partial_t|u^\varepsilon_{2}|^2\nonumber\\
&+&\frac{{\varepsilon a^\infty }}{2}
\int_0^t\int_{\Gamma^{wa}_\varepsilon}(\frac{1}{\delta}
|\partial_tu^\varepsilon_{2}|^2+\delta|\partial_tu^\varepsilon_{3}|^2)d\sigma_xd\tau\nonumber\\
&+&
\frac{\varepsilon b^\infty }{2}
\int_0^t\int_{\Gamma^{wa}_\varepsilon}
|\partial_tu^\varepsilon_{2}|^2d\sigma_xd\tau+\frac{\varepsilon B}{2}
\int_0^t\int_{\Gamma^{wa}_\varepsilon}
\partial_t|u^\varepsilon_{2}|^2d\sigma_xd\tau.\nonumber
\end{eqnarray}
Using (\ref{1a}), it yields
\begin{equation}\label{331}
\begin{aligned}
\frac{1}{2}\int_{\Omega^\varepsilon}|\partial_tu^\varepsilon_{2}|^2dx&+\left({d_2}-\frac{C^*A^\infty\varepsilon^2}{2}-\frac{C^*B^\infty\varepsilon^2}{2}-
\frac{C^*a^\infty\varepsilon^2}{2\delta}\right)
\int_0^t\int_{\Omega^\varepsilon}|\nabla \partial_t{u_{2}}^\varepsilon|^2dxd\tau
\\&\leq C_{13}+\frac{C^*a^\infty\delta}{2}\int_0^t\int_{\Omega^\varepsilon}(|\partial_tu^\varepsilon_{3}|^2+\varepsilon^2|\nabla\partial_tu^\varepsilon_{3}|^2)\\&+\left(  \frac{k_1^\infty}{2}+ \frac{K_1^\infty}{2}-k_2+ \frac{C^*A^\infty}{2}+\frac{C^*B^\infty\varepsilon^2}{2}+\frac{C^*a^\infty}{2\delta}\right)
\int_0^t\int_{\Omega^\varepsilon}|\partial_tu^\varepsilon_{2}|^2,
\end{aligned}
\end{equation}
where $C_{13}$ depends on the bounded terms of r.h.s of (\ref{aqw}).
Differentiating now ({\ref{def3}}) with respect to time and then
testing the result with $\partial_t(u^\varepsilon_{3}-u^D_{3})$, we
get
\vspace{-0.3in}
\begin{eqnarray*}
\frac{1}{2}\int_{\Omega^\varepsilon}|\partial_tu^\varepsilon_{3}|^2dx&+&{d_3}
\int_0^t\int_{\Omega^\varepsilon}|\nabla \partial_t{u_{3}}^\varepsilon|^2dxd\tau\\\nonumber&+&
\int_0^t\int_{\Omega^\varepsilon}\left(\frac{1}{2}(\partial_td_{3}|\nabla {u_{3}}^\varepsilon|^2-
(\partial_t\partial_td_{3})|\nabla{u_{3}}^\varepsilon|^2\right)\\\nonumber&\leq&
\frac{D_3}{2}\int_0^t\int_{\Omega^\varepsilon}|\nabla{u_{3}}^\varepsilon|^2dxd\tau+
\frac{d^\infty_3}{2}\int_0^t\int_{\Omega^\varepsilon}|\nabla\partial_t{u_{3}}^\varepsilon|^2dxd\tau\\\nonumber&
+&\frac{d^\infty_3+D_3}{2}\int_0^t\int_{\Omega^\varepsilon}|\nabla\partial_t{u_{3}}^D|^2dxd\tau+
\frac{\varepsilon A^\infty }{2}\int_0^t\int_{\Gamma^{wa}_\varepsilon}\partial_t|{u_{3}}^\varepsilon|^2dxd\tau\\\nonumber&+&
{\varepsilon a^\infty }{}\int_0^t\int_{\Gamma^{wa}_\varepsilon}|\partial_t{u_{3}}^\varepsilon|^2dxd\tau+
\frac{\varepsilon A^\infty }{2}\int_0^t\int_{\Gamma^{wa}_\varepsilon}(|{u_{3}}^\varepsilon|^2+|\partial_t{u_{3}}^D|^2)dxd\tau
\\\nonumber&+&\frac{\varepsilon a^\infty }{2}\int_0^t\int_{\Gamma^{wa}_\varepsilon}(|\partial_t{u_{3}}^\varepsilon|^2+|\partial_t{u_{3}}^D|^2)dxd\tau
\\\nonumber&+&\frac{\varepsilon B^\infty }{2}\int_0^t\int_{\Gamma^{wa}_\varepsilon}(|{u_{2}}^\varepsilon|^2+|\partial_t{u_{3}}^\varepsilon|^2+|{u_{2}}^\varepsilon|^2+|\partial_t{u_{3}}^D|^2)dxd\tau
\\\nonumber&+&\frac{\varepsilon b^\infty }{2}\int_0^t\int_{\Gamma^{wa}_\varepsilon}(\frac{1}{\delta}|\partial_t{u_{2}}^\varepsilon|^2+\delta|\partial_t{u_{3}}^\varepsilon|^2+
|\partial_t{u_{2}}^\varepsilon|^2+|\partial_t{u_{3}}^D|^2)dxd\tau.
\end{eqnarray*}
Using (\ref{1a}) to deal with the boundary terms, we obtain
\begin{eqnarray}\label{32}
\frac{1}{2}\int_{\Omega^\varepsilon}|\partial_tu^\varepsilon_{3}|^2dx
&+&\left({d_3}-\frac{d_3^\infty}{2}-\frac{C^*\varepsilon^2}{2}(3a^\infty
+B^\infty+b^\infty+a^\infty\delta)\right)
\int_0^t\int_{\Omega^\varepsilon}|\nabla \partial_t{u_{3}}^\varepsilon|^2dxd\tau
\nonumber
\end{eqnarray}
\begin{eqnarray}
&\leq&C_{14}+C_{15}\int_0^t\int_{\Omega^\varepsilon}|\partial_t{u_{3}}^\varepsilon|^2dxd\tau
\\&+&C^*b^\infty\int_0^t\int_{\Omega^\varepsilon}(|\partial_t{u_{3}}^\varepsilon|^2+\varepsilon^2|\nabla\partial_t{u_{3}}^\varepsilon|^2)dxd\tau
\end{eqnarray}
Adding (\ref{331}) and (\ref{32}) and using (\ref{29}) and (\ref{30}) to get the desired result.
\subsection{Extension step}
Since we deal here with an oscillating system posed in a perforated
domain, the natural next step is to extend all concentrations to the
whole $\Omega$. We do this by following a two-steps procedure: In
Step 1, we rely on the standard extension results indicated in
section \ref{extension-section} to extend all active concentrations
$u^\varepsilon_\ell$ ($\ell\in \{1,\dots,4\}$) to $\Omega$. In step
2, we unfold the ode for $u^\varepsilon_5$ such that the unfolded
concentration is defined on the fixed boundary $\Gamma$; see section
\ref{unfolding-section}.
\subsection{Extension lemmas}\label{extension-section}
Since all the concentrations are defined in $\Omega^\varepsilon$ and $\Omega^\varepsilon_1$, to get macroscopic equations we need to extend them into $\Omega$.
\begin{rem}
Take $\varphi^\varepsilon\in L^2(0,T;H^1(\Omega^\varepsilon))$. Note that since our microscopic geometry is sufficiently regular, we can speak in terms of extensions. Recall the linearity of the extension operator
$$\mathcal{P}^\varepsilon:L^2(0,T;H^1(\Omega^\varepsilon))\rightarrow L^2(0,T;H^1(\Omega))$$
defined by $\mathcal{P}^\varepsilon\varphi^\varepsilon=\tilde{\varphi}^\varepsilon$.
To keep notation simple, we denote the extension $\tilde{\varphi}^\varepsilon$ again by $\varphi^\varepsilon$.
\end{rem}

\begin{lem}{(Extension})\label{extension}
Consider the geometry described in Section \ref{geometry}. There
exists an extension $\tilde{u}^\varepsilon$ of $u^\varepsilon$ such
that
\begin{enumerate}
  \item  $\parallel\tilde{u}^\varepsilon\parallel_{L^2(Y)}\leq \hat{C} \parallel u^\varepsilon\parallel_{L^2(Y^w)},$  for $u^\varepsilon\in{L^2(Y^w)}$
  \item  $\parallel\nabla\tilde{u}^\varepsilon\parallel_{L^2(Y)}\leq \hat{C} \parallel \nabla u^\varepsilon\parallel_{L^2(Y^w)},$  for $ \nabla u^\varepsilon\in{L^2(Y^w)}$
  \item  $\parallel\tilde{u}^\varepsilon\parallel_{H^1(\Omega)}\leq \hat{C} \parallel u^\varepsilon\parallel_{H^1(\Omega^\varepsilon)}$,  for $u^\varepsilon\in H^1(\Omega^\varepsilon)$
\end{enumerate}
\end{lem}
\textit{Proof.}
For the proof of this Lemma, see Section 2 in \cite{paulin} or compare Lemma 5, p.214 in \cite{hornung}.

\begin{defn}(Two-scale convergence cf. \cite{Allaire,Nguestseng})\label{2 scale}
Let $\{u^\varepsilon\}$ be a sequence of functions in
$L^2((0,T)\times\Omega)$ ($\Omega$ being an open set of
$\mathbb{R}^N$) where $\varepsilon$ being a sequence of strictly
positive numbers that tends to zero. $\{u^\varepsilon\}$ is said to
two-scale converge to a unique function
 $u_0(t,x,y)\in L^2((0,T)\times\Omega\times Y)$ if and only if for any $\psi\in C^\infty_0((0,T) \times\Omega, C^\infty_\#(Y))$, we have
 \vspace{-0.3in}
\begin{eqnarray}\label{00}
lim_{\varepsilon \rightarrow 0}\int_0^T\int_{\Omega} u^\varepsilon \psi(t,x,\frac{x}{\varepsilon})dxdt = \int_{\Omega} \int_{Y} u_0(t,x,y) \psi(t,x,y) dydxdt.
\end{eqnarray}
We denote (\ref{00}) by $u^\varepsilon\stackrel{2}{\rightharpoonup}
u_0$.
\end{defn}
\begin{thm}
\begin{enumerate}
  \item[(i)] From each bounded sequence $\{u^\varepsilon\}$ in $L^2((0,T)\times\Omega)$, one can extract a
subsequence which two-scale converges to $u_0(t,x,y)\in L^2((0,T)\times\Omega\times Y)$.
  \item[(ii)] Let $\{u^\varepsilon\}$ be a bounded sequence in $H^1((0,T)\times\Omega)$, which converges weakly to a limit
function $u_0(t,x,y)\in H^1((0,T) \times \Omega \times Y)$. Then there exists $\tilde{u}\in L^2(\Omega;H^1_\#(Y)/\mathbb{R})$
 such that up to a subsequence $\{u^\varepsilon\}$ two-scale converges to $u_0(t,x,y)$ and
 $ \nabla u^\varepsilon\stackrel{2}{\rightharpoonup} \nabla_xu_0+ \nabla_y\tilde{u}.$
  \item[(iii)] Let $\{u^\varepsilon\}$ and $\{\varepsilon\nabla u^\varepsilon\}$ be bounded sequences in $L^2((0,T) \times \Omega)$,
   then there exists${u_0} \in L^2((0,T) \times \Omega;H^1_\#(Y))$ such that up to a subsequence $u^\varepsilon$ and $\varepsilon\nabla u^\varepsilon$ two-scale converge
to $u_0(t,x, y)$ and $\nabla_yu_0(t,x, y)$ respectively.
\end{enumerate}
\end{thm}

\begin{defn}(Two-scale convergence for $\varepsilon-$periodic hypersurfaces \cite{maria})\label{2 scale hyp}
A sequence of functions $\{u^\varepsilon\}$ in
$L^2((0,T)\times\Gamma_\varepsilon)$ is said to two-scale converge
to a limit $u_0\in L^2((0,T)\times\Omega \times \Gamma)$ if and only
if for any $\psi\in C^\infty_0((0,T) \times\Omega,
C^\infty_\#(\Gamma))$ we have
\vspace{-0.3in}
\begin{eqnarray*}
lim_{\varepsilon \rightarrow 0}\varepsilon\int_0^T\int_{\Gamma_\varepsilon} u^\varepsilon \psi(t,x,\frac{x}{\varepsilon})d\sigma_xdt =
 \int_{\Omega} \int_{\Gamma} u_0(t,x,y) \psi(t,x,y) d\sigma_ydxdt.
\end{eqnarray*}
\end{defn}
\begin{thm}\label{tt}
\begin{enumerate}
        \item[(i)] From each bounded sequence $\{u^\varepsilon\} \in L^2((0,T)\times\Gamma_\varepsilon)$, one can extract a
subsequence $u^\varepsilon$ which two-scale converges to a function $u_0 \in L^2((0,T)\times\Omega \times \Gamma)$.
        \item[(ii)] If a sequence of functions $\{u^\varepsilon\}$ is bounded in $L^\infty((0,T)\times\Gamma_\varepsilon)$,
        then $u^\varepsilon$ two-scale converges to a function $u_0 \in L^\infty((0,T)\times\Omega \times \Gamma)$.
      \end{enumerate}
\end{thm}
{\bf Proof}. For proof of (i), see \cite{maria} and the one for (ii), see \cite{MM}.

\begin{lem}\label{lem1}
Assume the hypotheses of Lemma \ref{lemma2} and Lemma
\ref{extension} to hold. The {\em a priori} estimates lead to the
following convergence results:
\begin{enumerate}
  \item[(a)] $u^\varepsilon_{i}\rightharpoonup u_{i}$ in $L^2(0,T;H^1(\Omega)$ for all $i\in\{1,2,3,4\}$,
  \item[(b)] $u^\varepsilon_{i}\stackrel{\ast}{\rightharpoonup} u_{i}$ in $L^\infty((0,T)\times \Omega)$,
  \item[(c)] $\partial_tu^\varepsilon_{i}\rightharpoonup \partial_tu_{i}$ in $L^2((0,T)\times \Omega)$,
  \item[(d)] $u^\varepsilon_{i}\rightarrow u_{i}$ in $L^2(0,T;H^\beta(\Omega))$ for $\frac{1}{2}<\beta<1$, also $\parallel u^\varepsilon_{i}- u_{i}\parallel_{L^2((0,T)\times \Gamma_\varepsilon)}\rightarrow 0$ as $\varepsilon\rightarrow 0$,
  \item[(e)] $u^\varepsilon_{i}\stackrel{2}{\rightharpoonup} u_{i}, \nabla u^\varepsilon_{i}\stackrel{2}{\rightharpoonup} \nabla_xu_{i}+ \nabla_yu_{i1}$,
      $u_{i1}\in L^2((0,T)\times \Omega;H^1_\#(Y)/\mathbb{R})$,
  \item[(f)] $ u^\varepsilon_{5}\stackrel{2}{\rightharpoonup} u_{5}$, and $u_{5}\in L^\infty((0,T)\times \Omega\times \Gamma^{sw})$,
  \item[(g)] $\partial_t u^\varepsilon_{5}\stackrel{2}{\rightharpoonup} \partial_tu_{5}$, and $u_{5}\in L^2((0,T)\times \Omega\times \Gamma^{sw})$.
\end{enumerate}
\end{lem}

{\em Proof.}
 (a) and (b) are obtained as a direct consequence of the fact that $u^\varepsilon_{i}$ is bounded
 in $L^2(0,T;H^1(\Omega))\cap L^\infty((0,T)\times \Omega)$;
up to a subsequence (still denoted by $u^\varepsilon_{i}$)
$u^\varepsilon_{i}$ converges weakly to $u_{i}$ in
$L^2(0,T;H^1(\Omega))\cap L^\infty((0,T)\times \Omega)$. A similar
argument gives (c). To get (d), we use the compact embedding
$H^{\beta'}(\Omega)\hookrightarrow H^\beta(\Omega)$, for $\beta\in
(\frac{1}{2},1)$ and $0<\beta<\beta'\leq1$ (since $\Omega$ has
Lipschitz boundary). We have
$$W:=\{u_i\in L^2(0,T;H^1(\Omega))\; \mbox{and}\;\partial_tu_i\in L^2((0,T)\times \Omega)\mbox{ for all}\;i\in\{1,2,3,4\}\}$$
For a fixed $\varepsilon$, $W$ is compactly embedded in
$L^2(0,T;H^\beta(\Omega))$ by the Lions-Aubin Lemma; cf. e.g. \cite{lions}.
Using the trace inequality (\ref{1b})
\vspace{-0.3in}
\begin{eqnarray*}
\parallel u^\varepsilon_i-u_i\parallel_{L^2((0,T)\times \Gamma_\varepsilon)}&\leq &{C^*_0}\parallel u^\varepsilon_i-u_i\parallel_{L^2(0,T;H^\beta(\Omega^\varepsilon))},\\
&\leq&C\parallel u^\varepsilon_i-u_i\parallel_{L^2(0,T;H^\beta(\Omega))},
\end{eqnarray*}
where $\parallel
u^\varepsilon_i-u_i\parallel_{L^2(0,T;H^\beta(\Omega))}\rightarrow
0$ as $\varepsilon\rightarrow 0.$ To investigate (e), (f) and (g),
we use the notion of two-scale convergence as indicated in
Definition \ref{2 scale} and \ref{2 scale hyp}. Since
$u^\varepsilon_i$ are bounded in $L^2(0,T;H^1(\Omega)$, up to a
subsequence $u^\varepsilon_i\stackrel{2}{\rightharpoonup}
 u_i$ in $L^2((0,T)\times\Omega\times Y)$, and $\nabla u^\varepsilon_i\stackrel{2}{\rightharpoonup} \nabla_xu_i+\nabla_y\tilde{u}_i$, $\tilde{u}_i\in
L^2((0,T)\times\Omega;H^1_\#(Y)/\mathbb{R})$. By Theorem \ref{tt},
$u^\varepsilon_{5}$ in $L^\infty((0,T)\times \Omega\times \Gamma)$
converges two-scale to $u_{5}$
 in the same space and $\partial_tu^\varepsilon_{5}$ converges two-scale to $\partial_tu_{5}$ in $L^2((0,T)\times \Omega\times \Gamma)$.
Due to the presence of the non-linear reaction rate on the interface
$\Gamma^{sw}_\varepsilon$, the convergences listed in Lemma
\ref{lem1} are still not sufficient to pass to the limit
$\varepsilon\rightarrow0$ in the microscopic model. To be more
precise, we can pass to $\varepsilon\rightarrow0$ in the pde's, but
not in the ode.

\subsection{Cell problems}

In order to be able to formulate the upscaled equations, we define
two classes of {\em cell problems} very much in the spirit of
\cite{Hornung2}. One class of problems will refer to the
water-filled parts of the pore, while the second class will refer to
the air-filled part of the pores.

\begin{defn}{(Cell problems)}\label{cell problems}
The \underline{cell problems in water-filled part} are given by
\begin{equation*}
\begin{aligned}
-\nabla_y.(D_\ell(t,y)\nabla_y\chi_i)&=\sum_{k=1}^3\partial_{y_k}{D_\ell}_{ki}(t,y), \;\mbox{in}\;\;Y^w,\\
-D_\ell(t,y)\frac{\partial \chi_i}{\partial n}&=\sum_{k=1}^3{D_\ell}_{ki}(t,y)n_k \;\;\mbox{on} \;\Gamma^{sw},\\
-D_\ell(t,y)\frac{\partial \chi_i}{\partial n}&=\sum_{k=1}^3{D_\ell}_{ki}(t,y)n_k \;\;\mbox{on} \;\Gamma^{wa},
\end{aligned}
\end{equation*}
for all $i, \ell\in\{1,2,4\}$ and $\chi_i$ are Y-periodic in y. The
\underline{cell problems in air-filled part} are given by
\begin{equation*}
\begin{aligned}
-\nabla_y.(D_3(t,y)\nabla_y\varsigma_i)&=\sum_{k=1}^3\partial_{y_k}{D_3}_{ki}(t,y), \;\mbox{in}\;\;Y^a,\\
-D_3(t,y)\frac{\partial \varsigma_i}{\partial n}&=\sum_{k=1}^3{D_3}_{ki}(t,y)n_k \;\;\mbox{on} \;\Gamma^{wa},\\
-D_3(t,y)\frac{\partial\varsigma_i}{\partial n}&=\sum_{k=1}^3{D_3}_{ki}(t,y)n_k \;\;\mbox{on} \;\partial Y^a-\Gamma^{wa},
\end{aligned}
\end{equation*}
for all $i\in\{1,2,3\}$ and $\varsigma_i$ are $Y$-periodic in $y$.
\end{defn}
\section{Two-scale limit equations}\label{two-scale}
\begin{thm}\label{thm}
The sequences of the solutions of the weak formulation
(\ref{def1})-(\ref{def5}) converges to the weak solution $u_i,
i\in\{1,2,3,4,5\}$ as $\varepsilon\rightarrow 0$ such that $u_i\in
H^1(0,T;L^2(\Omega))\cap L^2(0,T;H^1(\Omega))\cap
L^\infty((0,T)\times\Omega)$ and $u_5\in
H^1(0,T;L^2(\Omega\times\Gamma))\cap
L^\infty((0,T)\times\Omega\times\Gamma))$.
 The weak formulation of the two-scale limit equations is given by
 \vspace{-0.2in}
\begin{eqnarray}\label{macp}
\int_0^T\int_{\Omega}\partial_tu_{i}(t,x)\phi_i(t,x)dxdt&+&\int_0^T\int_{\Omega}\tilde{d}_i(t)\nabla u_{i}(t,x)
\nabla\phi_idxdt\\\nonumber
&=&\int_0^T\int_{\Omega}F_i(u)\phi_idxdt \mbox{ for all }i\in\{1,2,3,4\},
\end{eqnarray}
where
\vspace{-0.2in}
\begin{eqnarray*}
F_1(u)&:=&-\tilde{k}_1(t)u_1(t,x)+\tilde{k}_2(t)u_2(t,x)\\&&\;-\frac{1}{|Y|}\int_{\Gamma}k_3(t,y)R(u_1(t,x))Q(u_5(t,x,y))d\sigma_y,\\
F_2(u)&:=&\tilde{k}_1(t)u_1(t,x)-\tilde{k}_2(t)u_2(t,x)+\tilde{a}(t)u_3(t,x)-\tilde{b}(t)u_2(t,x),\\
F_3(u)&:=&-\tilde{a}(t)u_3(t,x)+\tilde{b}(t)u_2(t,x),\\
F_4(u)&:=&\tilde{k}_1(t)u_1(t,x),
\end{eqnarray*}
 with the initial values $ u_i(0,x)={u_i}_0(x)$ for $x\in \Omega$, and
\begin{equation}\label{maco}
\begin{aligned}
\int_0^T\int_{\Omega\times \Gamma}&\partial_t
u_5(t,x,y)\phi_5(t,x,y)dtdxd\sigma_y\\&=\int_0^T\int_{\Omega\times\Gamma}
k_3(t,y)R(u_1(t,x))Q(u_5(t,x,y))\phi_5(t,x,y)dt dx d\sigma_y,
\end{aligned}
\end{equation}
with $u_5(0,x,y)={u_5}_0(x,y)$ for $\; x\in \Omega,\;y\in\Gamma^{sw}$.
Also $\phi:=(\phi_1,\phi_2,\phi_3,\phi_4)\in[C^\infty((0,T)\times\Omega)]^4$, $\psi:=(\psi_1,\psi_2,\psi_3,\psi_4)\in[C^\infty((0,T)\times\Omega);C^\infty_\#(Y)]^4$,
\vspace{-0.2in}
\begin{eqnarray}
\tilde{k}_j(t)&:=&\frac{1}{|Y|}\int_Yk_j(t,y)dy,\;j\in\{1,2\},\\
\tilde{a}(t)&:=&\frac{1}{|Y|}\int_{\Gamma^{wa}}a(t,y)d\sigma_y,\\
\tilde{b}(t)&:=&\frac{1}{|Y|}\int_{\Gamma^{wa}}b(t,y)d\sigma_y,\\
{\tilde{d}_{\ell ij}}&:=&\sum_{k=1}^3\int_{Y}(d_{\ell ij}(t,y)+d_{\ell ik}(t,y)(\delta_{n\ell}\partial_{y_k}\chi_j+\delta_{3\ell}\partial_{y_k}\varsigma_j)dy,\\
&&\ell\in\{1,2,3,4\},\,n\in\{1,2,4\}\nonumber
\end{eqnarray}
 with $\chi_j, \varsigma_j$ being solutions of the cell
problems defined in Definition \ref{cell problems}, while $\delta$
denotes here the Kronecker's symbol.
\end{thm}
{\em Proof.} We apply two-scale convergence techniques together with
Lemma \ref{lem1} to get macroscopic equations. We take test
functions incorporating the following oscillating behavior
$\varphi_i(t,x)=\phi_i(t,x)+\varepsilon\psi_i(t,x,\frac{x}{\varepsilon}),
\phi_i\in C^\infty((0,T)\times\Omega), \phi_i\in
C^\infty((0,T)\times\Omega,; C^\infty_\#(Y)),i\in\{1,2,3,4\}$.
Applying two-scale convergence yields
\vspace{-0.2in}
\begin{eqnarray}\label{34}
|Y|\int_0^T\int_{\Omega}&&\partial_tu_{i}\phi_i(t,x)dxdt+\int_0^T\int_{\Omega}\int_{Y}d_i(t,y)(\nabla_xu_{i}(t,x)\nonumber
\\&+&\nabla_y\tilde{u}_{i}(t,x,y))
(\nabla_x\phi_i(t,x)+\nabla_y\psi_i(t,x,\frac{x}{\varepsilon}))dydxdt\nonumber\\&=&\int_0^T\int_{\Omega}f_i(u)\phi_i(t,x)dxdt.
\end{eqnarray}
\vspace{-0.2in}
\begin{eqnarray*}\label{}
\int_0^T\int_{\Omega}f_1(u)&&\phi_1(t,x)dxdt=-\lim_{\varepsilon\rightarrow 0}\int_0^T\int_{\Omega^\varepsilon}k_1^\varepsilon u^\varepsilon_1(\phi_1(t,x)
+\varepsilon\psi_1(t,x,\frac{x}{\varepsilon})) dx dt\\
&+&\lim_{\varepsilon\rightarrow 0}\int_0^T\int_{\Omega^\varepsilon}k_2^\varepsilon u^\varepsilon_2(\phi_1(t,x)
+\varepsilon\psi_1(t,x,\frac{x}{\varepsilon})) dx dt\\
&-&\lim_{\varepsilon\rightarrow 0} \varepsilon\int_0^T\int_{\Gamma^{sw}_\varepsilon}\eta(R(u^\varepsilon_1),Q(u^\varepsilon_5))(\phi_1(t,x)
+\varepsilon\psi_1(t,x,\frac{x}{\varepsilon}))d\sigma_xdt.
\end{eqnarray*}
Using Lemma \ref{lem1}, we have
\vspace{-0.3in}
\begin{eqnarray*}\label{}
\int_0^T\int_{\Omega}f_1(u)\phi_1(t,x)dxdt&=&-\int_0^T\int_{\Omega}\int_{Y}k_1(t,y)u_1(t,x)\phi_1(t,x) dy dx dt\\
&+&\int_0^T\int_{\Omega}\int_{Y}k_2(t,y)u_2(t,x)\phi_1(t,x) dy dx dt\\
&-&\lim_{\varepsilon\rightarrow 0} \varepsilon\int_0^T\int_{\Gamma^{sw}_\varepsilon}\partial_tu^\varepsilon_5(\phi_1(t,x)
+\varepsilon\psi_1(t,x,\frac{x}{\varepsilon}))d\sigma_xdt.
\end{eqnarray*}
\vspace{-0.2in}
\begin{eqnarray*}\label{}
\int_0^T\int_{\Omega}f_1(u)\phi_1(t,x)dxdt&=&-|Y|\int_0^T\int_{\Omega}\tilde{k}_1(t)u_1(t,x)\phi_1(t,x) dx dt\\
&+&|Y|\int_0^T\int_{\Omega}\tilde{k}_2(t)u_2(t,x)\phi_1(t,x) dx dt\\
&-&\int_0^T\int_{\Omega}\int_{\Gamma}\partial_tu_5\phi_1(t,x)d\sigma_ydxdt.
\end{eqnarray*}
\vspace{-0.2in}
\begin{eqnarray*}\label{}
\int_0^T\int_{\Omega}f_2(u)\phi_2(t,x)dxdt&=&\lim_{\varepsilon\rightarrow 0}\int_0^T\int_{\Omega^\varepsilon}k_1^\varepsilon u^\varepsilon_1(\phi_2(t,x)
+\varepsilon\psi_2(t,x,\frac{x}{\varepsilon})) dx dt\\
&-&\lim_{\varepsilon\rightarrow 0}\int_0^T\int_{\Omega^\varepsilon}k_2^\varepsilon u^\varepsilon_2(\phi_2(t,x)
+\varepsilon\psi_2(t,x,\frac{x}{\varepsilon})) dx dt\\
&+&\lim_{\varepsilon\rightarrow 0} \varepsilon\int_0^T\int_{\Gamma^{wa}_\varepsilon}a_\varepsilon u^\varepsilon_3(\phi_2(t,x)
+\varepsilon\psi_2(t,x,\frac{x}{\varepsilon}))d\sigma_xdt\\
&-&\lim_{\varepsilon\rightarrow 0} \varepsilon\int_0^T\int_{\Gamma^{wa}_\varepsilon}b_\varepsilon u^\varepsilon_2(\phi_2(t,x)
+\varepsilon\psi_2(t,x,\frac{x}{\varepsilon}))d\sigma_xdt.
\end{eqnarray*}
\vspace{-0.2in}
\begin{eqnarray*}\label{}
\int_0^T\int_{\Omega}f_2(u)\phi_2(t,x)dxdt&=&|Y|\int_0^T\int_{\Omega}\tilde{k}_1(t)u_1(t,x)\phi_2(t,x) dx dt\\\\
&-&|Y|\int_0^T\int_{\Omega}\tilde{k}_2(t)u_2(t,x)\phi_2(t,x) dx dt\\
&+&|Y|\int_0^T\int_{\Omega}\tilde{a}(t) u_3(t,x)\phi_2(t,x)dxdt\\
&-&|Y|\int_0^T\int_{\Omega}\tilde{b}(t) u_2(t,x)\phi_2(t,x)dxdt.
\end{eqnarray*}
\vspace{-0.2in}
We also have
\vspace{-0.2in}
\begin{eqnarray*}\label{}
\int_0^T\int_{\Omega}f_3(u)\phi_3(t,x)dxdt&=
&-|Y|\int_0^T\int_{\Omega}\tilde{a}(t) u_3(t,x)\phi_3(t,x)dxdt\\
&+&|Y|\int_0^T\int_{\Omega}\tilde{b}(t) u_2(t,x)\phi_3(t,x)dxdt
\end{eqnarray*}
and
\vspace{-0.2in}
\begin{eqnarray*}\label{}
\int_0^T\int_{\Omega}f_4(u)\phi_4(t,x)dxdt&=
&|Y|\int_0^T\int_{\Omega}\tilde{k}_1(t) u_1(t,x)\phi_4(t,x)dxdt.
\end{eqnarray*}
We set $\phi_i=0,i\in\{1,2,3,4\}$ in (\ref{34}) to calculate the
expression of the known function $\tilde{u}_1$ and obtain
$$\int_0^T\int_{\Omega}\int_{Y}d_i(t,y)(\nabla_xu_{i}(t,x)
+\nabla_y\tilde{u}_{i}(t,x,y))
\nabla_y\psi_i(t,x,\frac{x}{\varepsilon})dydxdt=0,\mbox{ forall
}\psi_i.$$ Since $\tilde{u}_1$ depends  linearly on $\nabla_xu_{1}$,
it can be defined as
$$\tilde{u}_i:=\sum_{j=1}^3\partial_{x_j}u_{i}(\delta_{in}\chi_j(t,y)+\delta_{3i}\varsigma_j(t,y))\;\mbox{for}\;n\in\{1,2,4\}$$
where the function $\chi_j, \varsigma_j$ are the unique solutions of
the cell problems defined in Definition \ref{cell problems}.
Setting $\psi_i=0$ in (\ref{34}), we get
\begin{equation*}
\begin{aligned}
\int_0^T\int_{\Omega}\int_{Y}&\sum_{j,k=1}^3{d_i}_{jk}(t,y)(\partial_{x_k} u_i(t,x)\\&+ \sum_{m=1}^3 (\delta_{in}\partial_{y_k}\chi_m+\delta_{3i}\partial_{y_k}\varsigma_m)\partial_{x_m}u_i(t,x))
\partial_{x_j}\phi_k(t,x)dydxdt\\&=|Y| \int_0^T\int_{\Omega}\sum_{j,k=1}^3 {\tilde{d}_{ijk}} \partial_{x_k} u_i(t,x) \partial_{x_j}\phi_i(t,x)dxdt.
\end{aligned}
\end{equation*}
Hence, the coefficients (entering the effective diffusion tensor) are given by
\begin{equation}
{\tilde{d}_{ijk}}:=\frac{1}{|Y|}\sum_{k=1}^3\int_{Y}(d_{\ell ij}(t,y)+d_{\ell ik}(t,y)(\delta_{in}\partial_{y_k}\chi_j+\delta_{3i}\partial_{y_k}\varsigma_j)dy.
\end{equation}
$i\in\{1,2,3,4\}$, $n\in\{1,2,4\}$ and $j,k\in\{1,2,3\}$.
\subsection{Passing to the limit $\varepsilon\to 0$ in
(\ref{def5})}\label{unfolding-section}
It is not yet possible to pass to the limit $\varepsilon\rightarrow 0$ with the convergence results stated in Lemma \ref{lem1}.
To overcome this difficulty, we use the notion of periodic unfolding. It si worth mentioning that there is an intimate  link between the two-scale convergence and weak convergence of the unfolded sequences;
 see \cite{Cioranescu2,MM}. The key idea is: Instead of getting strong convergence for $u^\varepsilon_5$, obtain strong convergence for the periodic unfolding of $u^\varepsilon_5$.
\begin{defn}
For $\varepsilon>0$, the boundary unfolding of a measurable function $\varphi$ posed on oscillating surface $\Gamma_\varepsilon$ is defined by
$$T^b_\varepsilon \varphi(x,y)=\varphi(\varepsilon y+\varepsilon k),\;\;y\in \Gamma, x\in\Omega$$
where $k:=[\frac{x}{\varepsilon}]$ denotes the unique integer combination $\Sigma_{j=1}^3k_je_j$ of the
periods such that $x-[\frac{x}{\varepsilon}]$ belongs to $Y$.{ Note that the oscillation due to the perforations are shifted into the second variable $y$ which belongs to fixed surface $\Gamma$.}
\end{defn}
\begin{lem}
If $u_\varepsilon$ converges two-scale to $u$ and $T_b^\varepsilon
u_\varepsilon$ converges weakly to $u^*$ in
$L^2((0,T)\times\Omega;L^2_{\#}(\Gamma))$, then $u=u^*$ a.e. in
$(0,T)\times \Omega\times\Gamma$.
\end{lem}
\vspace{-0.2in}
{\em Proof.} The proof details for this statement can be found in
Lemma 4.6 of \cite{MM}.
\begin{lem}\label{35}
If $\varphi\in L^2((0,T)\times\Gamma^\varepsilon)$, then the following identity holds
$$\frac{1}{|Y|}\|T^\varepsilon_b \varphi\|_{L^2((0,T)\times\Omega\times\Gamma)}=\varepsilon\|\varphi\|_{L^2((0,T)\times\Gamma^\varepsilon)}.$$
\end{lem}
\vspace{-0.2in}
{\em Proof.} Consider
\vspace{-0.2in}
\begin{equation*}
\begin{aligned}
\frac{1}{|Y|}|T^\varepsilon_b \varphi|^2_{L^2(\Omega\times\Gamma)}&=\frac{1}{|Y|}\int_{\Omega\times\Gamma}|T^\varepsilon_b \varphi|^2dxd\sigma_y
=\frac{1}{|Y|}\int_{\Omega\times\Gamma}T^\varepsilon_b \varphi^2dxd\sigma_y,\\
&=\frac{1}{|Y|}\Sigma^3_{k=1}\int_{\varepsilon(Y+k)}\int_{\Gamma}T^\varepsilon_b \varphi^2dxd\sigma_y
=\frac{1}{|Y|}\Sigma^3_{k=1}\int_{\varepsilon(Y+k)}dx\int_{\Gamma}\varphi^2d\sigma_y,\\
&=\Sigma^3_{k=1}\varepsilon^3\int_{\Gamma}\varphi^2d\sigma_y.
\end{aligned}
\end{equation*}
Changing variable $z=\varepsilon(y+k)$, where $k=[\frac{x}{\varepsilon}]$, we get
\vspace{-0.2in}
\begin{eqnarray*}
\frac{1}{|Y|}|T^\varepsilon_b \varphi|^2_{L^2(\Omega\times\Gamma)}
&=&\Sigma^3_{k=1}\varepsilon^3\int_{\Gamma}\varphi^2d\sigma_y=\Sigma^3_{k=1}\varepsilon\int_{\varepsilon(\Gamma+k)}
\varphi^2d\sigma_z=\varepsilon\int_{\Gamma^\varepsilon}\varphi^2d\sigma_z.
\end{eqnarray*}
This completes the proof of (\ref{35}).
\begin{lem}
If $\varphi\in L^2(\Omega)$, then $T^\varepsilon_b \varphi\rightarrow \varphi$ as $\varepsilon\rightarrow 0$ strongly in $L^2(\Omega\times\Gamma)$.
\end{lem}
{\em Proof.} See in \cite{Cioranescu1,Cioranescu3} for proof
details.

Using the boundary unfolding operator $T^\epsilon_b$, we unfold the ode (\ref{def5}).
Changing the variable, $x=\varepsilon y+\varepsilon k$ (for $x\in \Gamma^{sw}_\varepsilon$) to the fixed domain $(0,T)\times\Omega\times\Gamma$,
we have
\begin{eqnarray}\label{unfoldode}
\partial_tT^\epsilon_b u^\epsilon_5(t,x,y)=\eta(T^\epsilon_b u^{\epsilon }_1(t,x,y),T^{\epsilon}_b u^{\epsilon}_5(t,x,y)).
\end{eqnarray}
In the remainder of this section, we prove that $T^\varepsilon_b
u^\varepsilon_5$ converges strongly to $u_5$ in $L^2(\Omega
\times\Gamma)$. From the two-scale convergence of $u^\epsilon_5$, we
obtain weak convergence of $T^\epsilon u^\epsilon_5$ to $u_5$ in
$L^\infty((0,T)\times\Omega;L^2_{per}(\Gamma))$. We start with
showing that $\{T^\varepsilon_b u^\varepsilon_5\}$ is a Cauchy
sequence in $L^2(\Omega\times\Gamma)$. To this end, we choose
$m,n\in \mathbb{N}$ with $n>m$ arbitrary. Writing down
(\ref{unfoldode}) for the two different choices of $\varepsilon$
(i.e. $\varepsilon_i=\varepsilon_n$ and
$\varepsilon_i=\varepsilon_m$), we obtain after subtracting the
corresponding equations that
\vspace{-0.2in}
\begin{eqnarray}\label{37}
&\partial_t&\int_{\Omega\times{\Gamma}}|T^{\epsilon_n}_bu^{\epsilon_n}_5-T^{\epsilon_m}_bu^{\epsilon_m}_5|^2d\sigma_ydx\nonumber\\
&=&\int_{\Omega\times{\Gamma}}[k^\varepsilon_3 R(T^{\epsilon_n}_bu^{\epsilon_n}_1)Q(T^{\epsilon_n}_bu^{\epsilon_n}_5)-k^\varepsilon_3 R(T^{\epsilon_m}_bu^{\epsilon_m}_1)Q(T^{\epsilon_m}_b
u^{\epsilon_m}_5))\nonumber\\
&\times&(T^{\epsilon_n}_bu^{\epsilon_n}_5-T^{\epsilon_m}_bu^{\epsilon_m}_5)d\sigma_ydx,\nonumber\\
&\leq& k^\infty_3c_R(\frac{Q^\infty}{2}+c_Qsup_{\Omega\times{\Gamma}}|T^{\epsilon_n}_bu^{\epsilon_n}_1|)\int_{\Omega\times{\Gamma}}|T^{\epsilon_n}_bu^{\epsilon_n}_5-T^{\epsilon_m}_b
u^{\epsilon_m}_5|^2d\sigma_ydx\nonumber\\
&+&\frac{k^\infty_3c_RQ^\infty}{2}\int_{\Omega\times{\Gamma}}|T^{\epsilon_n}_bu^{\epsilon_n}_1-T^{\epsilon_m}_bu^{\epsilon_m}_1|^2d\sigma_ydx.
\end{eqnarray}
To get (\ref{37}), we have used the uniform boundedness of $T^{\epsilon_n}_bu^{\epsilon_n}_1$.
We consider now
\vspace{-0.3in}
\begin{eqnarray}\label{36}
\int_{\Omega\times{\Gamma}}&&|T^{\epsilon n}_bu^{\epsilon n}_1-T^{\epsilon n}_bu^{\epsilon m}_1|^2d\sigma_ydx\nonumber\\
&\leq&\int_{\Omega\times{\Gamma}}(|T^{\epsilon n}_bu^{\epsilon n}_1-T^{\epsilon m}_bu_1|^2+
|T^{\epsilon n}_bu_1-u_1|^2)d\sigma_ydx\nonumber\\&+&\int_{\Omega\times{\Gamma}}(|T^{\epsilon m}_bu_1-u_1|^2+
|T^{\epsilon m}_bu^{\varepsilon m}_1-T^{\epsilon m}_bu_1|^2)d\sigma_ydx.
\end{eqnarray}
Since $u_1$ is constant w.r.t. $y$, we have that $T^{\epsilon_m}_bu_1\rightarrow u_1$ strongly in $L^2((0,T)\times\Omega\times\Gamma)$ as $\varepsilon\rightarrow0$.
From Lemma \ref{35}, we conclude that
\vspace{-0.1in}
\begin{equation*}
\int_{\Omega\times{\Gamma}}|T^{\epsilon }_bu^{\epsilon }_1-T^{\epsilon }_bu_1|^2d\sigma_ydx
\leq\varepsilon \int_{{\Gamma}^\varepsilon}|u^{\epsilon }_1-u_1|^2d\sigma_ydx
\leq\varepsilon C.
\end{equation*}
(\ref{36}) turns out to be
\vspace{-0.3in}
\begin{eqnarray*}
\int_{\Omega\times{\Gamma}}&&|T^{\epsilon_n}_bu^{\epsilon_n}_1-T^{\epsilon_n}_bu^{\epsilon_m}_1|^2d\sigma_ydxdt
\leq C(\varepsilon_n+\varepsilon_m),
\end{eqnarray*}
while (\ref{37}) becomes
\begin{equation*}
\partial_t\int_{\Omega\times{\Gamma}}|T^{\epsilon n}_bu^{\epsilon n}_5-T^{\epsilon m}_bu^{\epsilon m}_5|^2d\sigma_ydx
\leq C_{15}\int_{\Omega\times{\Gamma}}|T^{\epsilon_n}_bu^{\epsilon_n}_5-T^{\epsilon_m}_bu^{\epsilon_m}_5|^2d\sigma_ydx+\frac{C_{16}}{n},
\end{equation*}
where $C_{15}:=k^\infty_3c_R(\frac{Q^\infty}{2}+c_Qsup_{\Omega\times{\Gamma}}|T^{\epsilon_n}_bu^{\epsilon_n}_1|)$ and $C_{16}:=\frac{k^\infty_3c_RQ^\infty}{2}C$.
The Gronwall's inequality gives
\vspace{-0.3in}
\begin{eqnarray}\label{666}
\parallel T^{\epsilon n}_bu^{\epsilon n}_5-T^{\epsilon m}_bu^{\epsilon m}_5\parallel_{L^2(\Omega\times\Gamma)}\leq \frac{C_{16}}{n}.
\end{eqnarray}
By (\ref{666}), $\{T^{\epsilon}_bu^{\epsilon}_5\}$ is a Cauchy sequence.
Now, we take the two-scale limit in the ode (\ref{unfoldode}) to get
\vspace{-0.3in}
\begin{eqnarray*}\label{38}
\lim_{\varepsilon\rightarrow0}\varepsilon\int^T_0\int_{\Gamma^{sw}_\varepsilon}
\partial_t T^{\epsilon }_bu^{\varepsilon}_5\phi_1(t,x,\frac{x}{\varepsilon})d\sigma_xdt
=\lim_{\varepsilon\rightarrow0}\varepsilon\int_0^T\int_{\Gamma^{sw}_\varepsilon} \eta(T^{\epsilon }_bu^{\varepsilon}_1,T^{\epsilon }_bu^{\varepsilon}_5)\phi_1(t,x,\frac{x}{\varepsilon})d\sigma_xdt\nonumber.
\end{eqnarray*}
Consequently, we have
\vspace{-0.2in}
\begin{eqnarray}\label{66}
\int^T_0\int_{\Omega \times\Gamma^{sw}}&&
\partial_tu_5\phi_5(t,x,y)dxd\sigma_ydt\nonumber\\&=&\lim_{\varepsilon\rightarrow0}\varepsilon\int^T_0
\int_{\Gamma^{sw}_{\varepsilon}}T^{\epsilon }_b k^{\varepsilon}_3 R(T^{\epsilon }_bu^{\varepsilon}_1)Q(u^{\varepsilon}_5)\phi_5(t,x,\frac{x}{\varepsilon})
d\sigma_xdt,\nonumber\\
&=&\lim_{\varepsilon\rightarrow0}\varepsilon\int^T_0\int_{\Gamma^{sw}_{\varepsilon}}T^{\epsilon }_b k^{\varepsilon}_3 R(T^{\epsilon }_bu^{\varepsilon}_1)
Q(u_5)\phi_5(t,x,\frac{x}{\varepsilon})d\sigma_xdt\nonumber
\\&+&\lim_{\varepsilon\rightarrow0}\varepsilon\int^T_0\int_{\Gamma^{sw}_{\varepsilon}}T^{\epsilon }_b k^{\varepsilon}_3 R(T^{\epsilon }_bu^{\varepsilon}_1)
(Q(T^{\epsilon }_bu^{\varepsilon}_5)-Q(u_5))\phi_5(t,x,\frac{x}{\varepsilon})d\sigma_xdt.\nonumber\\
\end{eqnarray}
By (A2) and the strong convergence of $u^{\varepsilon}_1$, the first
term on the right hand side of (\ref{66}) converges two-scale to
$$\int^T_0\int_\Omega\int_{\Gamma^{sw}} k_3(t,y)R(u_1)Q(u_5)\phi_5(t,x,y)d\sigma_ydxdt,$$
while the second integral of (\ref{66})
\vspace{-0.2in}
\begin{eqnarray*}
&\varepsilon&\int^T_0\int_{\Gamma^{sw}_{\varepsilon}} T^{\epsilon }_bk^{\varepsilon}_3 R(T^{\epsilon }_bu^{\varepsilon}_1)
(Q(T^{\epsilon }_bu^{\varepsilon}_5)-Q(u_5))\phi_5(t,x,\frac{x}{\varepsilon})d\sigma_xdt\\
&\leq &\varepsilon\left(\int^T_0\int_{\Gamma^{sw}_{\varepsilon}} |T^{\epsilon }_bk^{\varepsilon}_3 R(T^{\epsilon }_bu^{\varepsilon}_1)
\phi_5(t,x,\frac{x}{\varepsilon})|^2d\sigma_xdt\right)^{\frac{1}{2}}\cdot \nonumber\\
&\cdot&\left(\int^T_0\int_{\Gamma^{sw}_{\varepsilon}}
|Q(T^{\epsilon }_bu^{\varepsilon}_5)-Q(u_5)|^2d\sigma_xdt\right)^{\frac{1}{2}},\\
&\rightarrow&0 \;\mbox{as}\; \varepsilon\rightarrow 0.
\end{eqnarray*}
At this point, we have used again (A2) in combination with the
strong convergence of $ T^{\epsilon }_bu^\varepsilon_5$. So, as
result of passing to the limit $\varepsilon\rightarrow0$ in
(\ref{def5}) we get (\ref{maco}).

It is worth noting that the weak
solution to the two-scale model inherits a.e. the positivity and
boundedness properties
 from the corresponding properties of the weak solution of the microscopic model.  Now, it only remains to ensure the uniqueness of weak solutions to the upscaled model.
\begin{lem}(Uniqueness of solutions of (\ref{macp})-(\ref{maco})
Assume (A1)-(A6). There exists at most one weak solution to the two-scale limit problem (\ref{macp}) and (\ref{maco}).
\end{lem}
\vspace{-0.2in}
{\em Proof.} Suppose there are two weak solutions to the two-scale
limit problem $(u^j_1,u^j_2,u^j_3,u^j_4,u^j_5)$ with $j\in\{1,2\}$.
We denote $u_\ell=u^1_\ell-u^2_\ell,\;\ell\in\{1,2,3,4\}$ and choose
as test function $\phi_\ell=u_\ell$. After straightforward
calculations, we have from (\ref{maco})
\vspace{-0.2in}
\begin{eqnarray}\label{un11}
|u^1_5-u^2_5| \leq C \int_0^t |u^1_1-u^2_1| d\tau.
\end{eqnarray}
Take $\phi_1=u_1$ in (\ref{macp}) to obtain
\vspace{-0.2in}
\begin{eqnarray}\label{un2}
\frac{1}{2} \int^t_0 \int_\Omega \partial_t |u_1|^2 dx dt& +& \tilde{d}_1 \int^t_0 \int_\Omega |\nabla u_1|^2 dx dt \nonumber\\
&\leq& -\tilde{k}_1 \int^t_0 \int_\Omega \partial_t |u_1|^2 dx dt +\frac{\hat{k}_2^\infty}{2}\int^t_0
\int_\Omega  (|u_1|^2 + |u_2|^2) dx dt \nonumber\\&
+&\frac{k_3^\infty}{|Y|} c_R c_Q M_1 \int^t_0 \int_{\Omega \times \Gamma^{sw}} (u^1_5-u^2_5)u_1 dx d\sigma_y dt \nonumber\\&
+& \frac{k_3^\infty}{|Y|} c_R Q^\infty \int^t_0 \int_{\Omega \times \Gamma^{sw}} |u_1|^2 dx d\sigma_y dt.
\end{eqnarray}
Using (\ref{un11}) together with the trace inequality for fixed domains, see section 5.5 Theorem 1 in \cite{Evans} and
 also the fact that $u_1$ is independent of y in (\ref{un2}), we get
 \vspace{-0.2in}
\begin{eqnarray*}
 \int^t_0 \int_\Omega \partial_t |u_1|^2 dx dt& +& (2\tilde{d}_1 - {k_3^\infty c_R}{} C^*
 (\delta M_1+Q^\infty))\int^t_0 \int_\Omega |\nabla u_1|^2 dx d\tau \nonumber\\
 &+&2\tilde{k}_1 \int^t_0 \int_\Omega \partial_t |u_1|^2 dx d\tau \nonumber\\ & \leq &
 ({\tilde{k}_2^\infty}+ {k_3^\infty c_R}{} C^* (\delta M_1+Q^\infty))
\int^t_0 \int_\Omega (|u_1|^2 + |u_2|^2) dx d\tau \nonumber\\&
+&{k_3^\infty}{\delta } c_R c_QM_1C^* \int^t_0 \int_{\Omega } \int^\tau_0(|u_1|^2 + |\nabla u_1|^2 ) ds dxd\tau .
\end{eqnarray*}
 For suitable choice of $\delta \in ]0, \frac{2d_1-k_3^\infty c_RC^*Q^\infty}{k_3^\infty c_RC^*M_1}[$, we have
 \vspace{-0.2in}
 \begin{eqnarray}\label{un3}
 \int^T_0 \int_\Omega \partial_t |u_1|^2 dx dt& +& \tilde{d}_1 \int^T_0 \int_\Omega |\nabla u_1|^2 dx d\tau
 +2\tilde{k}_1 \int^T_0 \int_\Omega \partial_t |u_1|^2 dx d\tau \nonumber\\ & \leq &
 ({\tilde{k}_2^\infty}+ k_3^\infty c_R C^* (\delta M_1+Q^\infty))
\int^T_0 \int_\Omega (|u_1|^2 + |u_2|^2) dx d\tau \nonumber\\&
+&\frac{k_3^\infty}{\delta} c_R c_Q M_1C^* \int^T_0 \int_{\Omega } \int^\tau_0(|u_1|^2 + |\nabla u_1|^2 ) ds dx d\tau .
\end{eqnarray}
Take $\phi_2=u_2$ in (\ref{macp}), we get
\vspace{-0.2in}
\begin{eqnarray*}\label{}
\frac{1}{2} \int^t_0 \int_\Omega \partial_t |u_2|^2 dx d\tau& +& \tilde{d}_2 \int^t_0 \int_\Omega |\nabla u_2|^2 dx d\tau \nonumber\\
&\leq& -\tilde{k}_2 \int^t_0 \int_\Omega \partial_t |u_2|^2 dx d\tau +\frac{\tilde{k}_1^\infty}{2}\int^t_0 \int_\Omega
(|u_1|^2 + |u_2|^2) dx d\tau \nonumber\\&
+&{\tilde{a}^\infty} \int^t_0 \int_{\Omega } u_2u_3 dx dt - {\tilde{b}} \int^T_0 \int_{\Omega } |u_2|^2 dx d\tau.
\end{eqnarray*}
\vspace{-0.3in}
\begin{eqnarray}\label{un4}
 \int^t_0 \int_\Omega \partial_t |u_2|^2 dx d\tau& +& \tilde{d}_2 \int^t_0 \int_\Omega |\nabla u_2|^2 dx d\tau \nonumber\\
&\leq& ({\tilde{k}_1^\infty} + \tilde{a}^\infty) \int^t_0 \int_\Omega (|u_1|^2 + |u_2|^2 +  |u_3|^2 ) dx d\tau .
\end{eqnarray}
Similarly, we obtain from (\ref{macp})
\vspace{-0.2in}
\begin{eqnarray}\label{un5}
 \int^t_0 \int_\Omega \partial_t |u_3|^2 dx d\tau& +& \tilde{d}_3 \int^t_0 \int_\Omega |\nabla u_3|^2 dx d\tau \leq
 \tilde{b}^\infty \int^t_0 \int_\Omega ( |u_2|^2 +  |u_3|^2 ) dx d\tau ,
\end{eqnarray}
\vspace{-0.3in}
\begin{eqnarray}\label{un6}
 \int^t_0 \int_\Omega \partial_t |u_4|^2 dx d\tau& +& \tilde{d}_4 \int^t_0 \int_\Omega |\nabla u_4|^2 dx d\tau \leq
   \tilde{k}^\infty_1 \int^T_0 \int_\Omega ( |u_1|^2 +  |u_3|^2 ) dx d\tau .
\end{eqnarray}
Adding side by side (\ref{un3})-(\ref{un6}) and applying Gronwall's inequality to the corresponding result, we receive
\vspace{-0.2in}
\begin{eqnarray}\label{a11}
 \Sigma_{i=1}^4\int_\Omega  |u_i|^2 dx + \hat{d}\Sigma_{i=1}^4\int^t_0 \int_\Omega |\nabla u_i|^2 dx d\tau +
   \hat{d} \int^t_0 \int_\Omega  |u_1|^2  dx d\tau \leq 0.
\end{eqnarray}
In (\ref{a11}), we have $\hat{d}:=min\{\tilde{d}_1,\tilde{d}_2,\tilde{d}_3,\tilde{d}_4,\tilde{k}_1\}>0$. Taking in (\ref{b11}) supremum over $(0,T)$, we obtain
\vspace{-0.2in}
\begin{eqnarray}\label{b11}
 \Sigma_{i=1}^4\int_\Omega  |u_i|^2 dx + \hat{d}\Sigma_{i=1}^4\int^T_0 \int_\Omega |\nabla u_i|^2 dx d\tau  \leq 0,
\end{eqnarray}
which concludes the proof of the Lemma.
\begin{lem}\label{st}(Strong formulation of the two-scale limit equations)
Assume the hypothesis of Lemma \ref{lem1} to hold. Then the strong
formulation of the two-scale limit equations (for all $t\in(0,T)$)
reads
\begin{equation}
\begin{aligned}
\partial_tu_{1}(t,x)&+\nabla\cdot(-\tilde{d}_1\nabla u_{1}(t,x))\\
&= -\tilde{k}_1(t)u_1(t,x)+\tilde{k}_2(t)u_2(t,x)\\&-\frac{1}{|Y|}\int_{\Gamma^{sw}}k_3(t,y)R(u_1(t,x))Q(u_5(t,x,y))d\sigma_y,\;x\in\Omega\\
u_1(0,x)&=u_{10}(x), \;x\in\bar\Omega,\\
n\cdot(-\tilde{d}_1\nabla u_{1}(t,x))&=0, \;x\in\partial\Omega
\end{aligned}
\end{equation}
\begin{equation}
\begin{aligned}
\partial_tu_{2}(t,x)&+\nabla\cdot(-\tilde{d}_2\nabla u_{2}(t,x))
= \tilde{k}_1(t)u_1(t,x)-\tilde{k}_2(t)u_2(t,x)\\&+\tilde{a}(t)u_3(t,x)-\tilde{b}(t)u_2(t,x),\;x\in\Omega,\\
u_2(0,x)&=u_{20}(x), \;x\in\bar\Omega,\\
n\cdot(-\tilde{d}_2\nabla u_{2}(t,x))&=0, \;x\in\partial\Omega,
\end{aligned}
\end{equation}
\begin{equation}
\begin{aligned}
\partial_tu_{3}(t,x)&+\nabla\cdot(-\tilde{d}_3\nabla u_{3}(t,x))
= -\tilde{a}(t)u_3(t,x)+\tilde{b}(t)u_2(t,x),\;x\in\Omega,\\
u_3(0,x)&=u_{30}(x), \;x\in\bar\Omega,\\
u_3(t,x)&=u^D_{3}(x), \;x\in\Gamma^{D},\\
n\cdot(-\tilde{d}_3\nabla u_{3}(t,x))&=0, \;x\in\Gamma^{N}
\end{aligned}
\end{equation}
\begin{equation}
\begin{aligned}
\partial_tu_{4}(t,x)&+\nabla\cdot(-\tilde{d}_4\nabla u_{4}(t,x))
= \tilde{k}_1(t)u_1(t,x),\;x\in\Omega,\\
u_4(0,x)&=u_{40}(x), \;x\in\bar \Omega,\\
n\cdot(-\tilde{d}_4\nabla u_{4}(t,x))&=0,\;x\in\partial\Omega
\end{aligned}
\end{equation}
\begin{equation}
\begin{aligned}
\partial_tu_5(t,x,y)&=k_3(t,y)R(u_1(t,x))Q(u_5(t,x,y)),\;x\in\Omega,y\in\Gamma^{sw},\\
u_5(0,x,y)&=u_{50}(x,y)\;x\in\bar\Omega,y\in\Gamma^{sw},
\end{aligned}
\end{equation}
where $\tilde{d}_i, i\in\{1,2,3,4\}$ and $\tilde{k}_j, j\in\{1,2\}$
are defined in Theorem \ref{thm}.
\end{lem}

\vspace{-0.3in}
\subsection*{Acknowledgements}

\vspace{-0.3in} We would like to thank M. Ptashnyk (RWTH Aachen) and
M. A. Peletier (TU Eindhoven) for fruitful discussions on this
subject.

\bibliography{tfam}
\bibliographystyle{elsart-num}

\end{document}